\documentclass[11 pt]{article}
\usepackage[top=2cm, bottom=2cm, left=25mm, right=25mm]{geometry}

\usepackage[utf8]{inputenc}
\usepackage{amsmath}
\usepackage{parskip}
\usepackage{amssymb}
\usepackage{color}
\usepackage[all]{xy}  
\usepackage{tikz-cd}
\usepackage{float}
\usepackage{setspace}
\usepackage{hyperref}
\newtheorem{theorem}{Theorem}

\newtheorem{proposition}[theorem]{Proposition}

\newtheorem{example}{Example}
\newtheorem{definition}{Definition}

\usepackage{caption}
\usepackage{lineno}

\usepackage{algorithmicx}
\usepackage{algorithm}
\usepackage{algpseudocode}

\usepackage[acronym]{glossaries}
\newacronym{RKMK}{RKMK}{Método de Runge-Kutta-Munthe-Kaas}

\usepackage{hyperref}
\hypersetup{colorlinks = true}

\usepackage{tikz-cd}
\usepackage{pgfplots}
\usepgfplotslibrary{groupplots}
\pgfplotsset{compat = 1.16}
\usetikzlibrary{decorations.markings}

\author{L. Blanco$^{*}$, F. Jiménez$^{*}$, J. de Lucas$^\dagger$, C. Sardón$^{*}$}

\begin{document}
\centerline{\Large Geometric numerical  methods  for Lie systems}
\vskip .3cm
\centerline{\Large and their application in optimal control}
\vskip .5cm

\centerline{L. Blanco$^{*}$, F. Jiménez$^{*}$, J. de Lucas$^\dagger$, C. Sardón$^{*}$}
\vskip .2cm

\centerline{$^{\dagger}$Department of Mathematical Methods in Physics, University of Warsaw,}
	\medskip
	\centerline{ul. Pasteura 5, 02-093, Warsaw, Poland.}
	\medskip
	\centerline{$^{*}$ Department of Applied Mathematics, Universidad Polit\'ecnica de Madrid (UPM)}
	\medskip
	\centerline{c. José Gutiérrez Abascal 2, 28006, Madrid.}
	\medskip

\begin{abstract}
A Lie system is a nonautonomous system of first-order ordinary differential equations whose general solution can be written via an autonomous function, the so-called (nonlinear) superposition rule of a finite number of particular solutions and some parameters to be related to initial conditions. This superposition rule can be obtained using the geometric features of the Lie system, its symmetries, and the symmetric properties of certain morphisms involved. Even if a superposition rule for a Lie system is known, the explicit analytic expression of its solutions frequently is not. This is why this article focuses on a novel geometric attempt to integrate Lie systems analytically and numerically. We focus on two families of methods based on Magnus expansions and on  Runge-Kutta-Munthe-Kaas methods, which are here adapted, in a geometric manner, to Lie systems. To illustrate the accuracy of our techniques we analyse Lie systems related to Lie groups of the form SL$(n,\mathbb{R})$, which play a very relevant role in mechanics. In particular, we depict an optimal control problem for a vehicle with quadratic cost function. Particular numerical solutions of the studied examples are given.

\end{abstract}

{\it MSC 2020 classes: 34A26; 53A70; (primary) 37M15; 49M25 (secondary)}

\section{Introduction}

The analytic integration of differential equations can be achieved in many relevant occasions, but it is not the usual case. Sometimes the geometric and symmetry properties of a Lie system are not enough to completely integrate the system, and this is why numerical methods are so important to study solutions of  differential equations. In particular, this paper devises geometric numerical methods adapted to a particular class of nonautonomous first-order systems of ordinary differential equations (ODEs): the so-called {\it Lie systems} \cite{carinena00,araujo20,win82}. 

A Lie system is a nonautonomous first-order system of ODEs that admits a general solution in terms of an autonomous function, the so-called {\it superposition rule}, a family of generic particular solutions and certain constants of integration related to the initial conditions \cite{carinena01, lucas11, sardon15}.  It is worth noting that a superposition rule for a Lie system may be explicitly known even when the explicit expression of its analytic solution is not \cite{lucas11}. Although obtaining a superposition rule reduces the integration of Lie systems to obtaining some particular solutions, 
 such particular solutions are not easy to describe explicitly \cite{lucas11,win82}. This is why we consider that geometric numerical methods for Lie systems should be developed. One could find an extensive list of works devoted to numerical algorithms in geometric mechanics \cite{CortesMartinez, iserles05,MMM, MaWe,McLac,SS} and references therein, but, as far as we know, just a few methods have been specifically designed for Lie systems \cite{piet12,RW84}. This manuscript, therefore, provides a novel application of geometric analytical and numerical methods to Lie systems, leading to some interesting consequences.

Our interest in Lie systems is two-fold. On the one hand, it is rooted in their geometric background. Long story short,  the origin of Lie systems goes back to the XIX century, when Sophus Lie proved that a nonautonomous system of ODEs of first-order admits a superposition rule if and only if it describes the integral curves of a $t$-dependent vector field defined taking values in a finite dimensional Lie algebra of vector fields, known as a {\it Vessiot-Guldberg Lie algebra} (VG henceforth) of the Lie system. The symmetries of a Lie system are direct correlation with the underlying VG Lie algebra. The theory of Lie systems has been widely studied in the last two decades and its research involves projective foliations, generalized distributions, Lie group theory, Poisson coalgebras, etc. (see \cite{carinena00,carinena07,araujo20} and references therein). In particular, the coalgebra method is based in symmetric properties of certain operators that allow us to obtain superposition rules with the aid of a finite-dimensional Poisson algebra of functions. 
On the other hand, Lie systems have many remarkable applications in many relevant scientific fields (see \cite{araujo20} and references therein). For instance, Lie systems are used in the study of the integrability of Riccati equations  \cite{carinena99}, quantum mechanics \cite{angelo05}, stochastic mechanics \cite{lazaro09}, superequations \cite{hussin90}, in biology and cosmology \cite{araujo20}.  Recently, the theory of Lie systems has been generalized to higher-order ordinary differential equations, such as higher-order Riccati equations \cite{carinena11}, second- and third-order Kummer-Schwarz equations \cite{lucas11},  Milne-Pinney equations \cite{carinena09}, among others. Additionally, the theory of Lie systems is also extensible to systems of partial differential equations \cite{carinena07, odzijewicz00}.

In the past few decades, discrete methods have made big progress in faithfully describing reality. For instance, the interest of numerical analysis in the research on Lie systems was already stressed by Winternitz \cite{win82}, who remarked that superposition rules allow us to study all solutions of a Lie system from the knowledge of some of them, which can be derived numerically. This is why the discretization of Lie systems and their numerical integration has caught our attention.
Since Lie systems are geometrically described in terms of an underlying VG Lie algebra, this allowed for solving a Lie system by studying a Lie system of a specific type, a so-called {\it automorphic Lie system} \cite{araujo20}, on a Lie group associated to the VG Lie algebra. Two automorphic Lie systems are two Lie systems that are equivalent under  automorphic transformations. In this way, automorphic Lie systems can be claimed to be symmetric Lie systems. One can then propose a numerical method for the automorphic Lie system, giving rise to numerical methods for a plethora of Lie systems that are related to the initial through an automorphic map that preserves the properties of the Lie group, aka. symmetry group transformation. Our perspective here on numerical methods specifically designed for Lie systems proposes numerical schemes on the Lie group.
There already exist some numerical methods designed to work on Lie groups, but our aim is to adapt them for Lie systems. In particular, we will focus on two classes of methods: the so-called Magnus methods \cite{iserles05, iserles99, Zanna} and Runge-Kutta-Munthe-Kaas (RKMK) \cite{munthe98, munthe99}, the latter being based on the classical Runge-Kutta (RK) schemes.

Summarizing, this manuscript presents a novel procedure for the integration of Lie systems by applying geometric numerical methods on one of its associated {\it automorphic Lie systems}, which is defined on a Lie group (we may refer to it as a VG Lie group). We aim at providing a quantitative and qualitative analysis of our numerical methods on the Lie group and compare them with the results obtained from numerical integration of the system of ODEs that defines the Lie system. This would resolve at the same time all Lie systems that are related to the same automorphic Lie system, i.e. those Lie systems that have isomorphic VG Lie algebras and that are determined by an equivalent curve within them (see \cite{carinena00} for details). We apply our numerical methods to automorphic Lie systems defined on Lie groups $\mathrm{SL}(n, \mathbb{R}) $, which appear in many physical applications (cf. \cite{LG17}). We are particularly interested in control theory, which involves matrix Riccati equations \cite{Ku73}.
We depict an application of matrix Riccati equations in optimal control with quadratic cost functions and solve it numerically with our adapted Magnus and RKMK methods.

 The structure of the paper goes as follows. Section 2 surveys the basic theory of Lie systems and develops their analytical resolution constructed upon the geometric structure they are built on; this analytical solution is enclosed in the procedure that is summarized in Procedure \ref{Metodo7Pasos}. Meanwhile, Section 3 is concerned with the novel discretization we are proposing for Lie systems, enclosed in Definition \ref{MetodosLS}. An application of our methods to SL$(2,\mathbb{R})$ and SL$(3,\mathbb{R})$  is provided in Section 4. Meanwhile, an optimal control problem for a vehicle with a quadratic cost function is presented in Section 5, and resolved using the novel analytical techniques we are delivering. 
\section{Geometric fundamentals and Lie systems} \label{fundamentos}

This section establishes the notation and geometric fundamentals on Lie systems and related concepts that we will be using throughout the manuscript. Unless otherwise stated, we hereafter assume all structures to be smooth, real  and globally defined. This will simplify the presentation, while stressing its main points. From now on, $\mathbb{K}$ stands for a field to be $\mathbb{R}$ or $\mathbb{C}$. 

\subsection{Geometric fundamentals}

A key concept in the theory of Lie systems is that of $t$-dependent vector fields. Let us describe  this geometric concept. Consider an $n$-dimensional manifold $N$ and its natural tangent bundle projection $\pi_N:TN\rightarrow N$. Let us define the projection
$ \pi_2 : (t, x) \in \mathbb{R} \times N \mapsto x \in N $, where $ t $ is the natural coordinate system on $ \mathbb{R} $. A {\it $t$-dependent vector field} on $ N $ is a map  $ X : (t, x) \in \mathbb{R} \times N \mapsto X(t, x) \in TN $ so that the following  diagram becomes commutative
\begin{equation*}
    \begin{tikzcd}
    \mathbb{R}\times N \arrow{r}{X} \arrow{dr}{\pi_2} & TN \arrow{d}{\pi_N} \\
                                                      & N
    \end{tikzcd}
\end{equation*}
i.e., $ \pi_N \circ X = \pi_2 $. 
In other words, a $t$-dependent vector field $X$ on $N$ amounts to a $t$-parametrized family of standard vector fields $\{X_t:x\in N\mapsto X(t,x)\in TN\}_{t\in \mathbb{R}}$ on $N$ (see \cite{lucas11} for details). We write $\mathfrak{X}_t(N)$ for the space of $t$-dependent vector fields on $N$, while $\mathfrak{X}(N)$ stands for the space of vector fields on $N$. 

An {\it integral curve} of a $t$-dependent vector field $X$ on $N$ is a curve $\gamma:\mathbb{R}\rightarrow  N$ of the form  $\gamma=\pi_2\circ \widehat{\gamma}$, where $\widehat{\gamma}:\mathbb{R}\rightarrow \mathbb{R}\times N$ is an integral curve of the so-called {\it autonomization} of $ X $, namely the vector field $\overline{X}=\partial/\partial t+ X$ on $\mathbb{R}\times N$, that is also a section of the natural projection $\pi:(t,x)\in \mathbb{R}\times N\mapsto t \in \mathbb{R}$. More precisely, if 
$X = \sum_{i = 1}^n \eta_i(t, x) \frac{\partial}{\partial x_i}$ on a local coordinate system $\{x_1,\ldots,x_n\}$ on $N$, 
then
$$\overline{X} = \frac{\partial}{\partial t} + \sum_{i = 1}^n \eta_i(t, x) \frac{\partial}{\partial x_i}
$$
and 
$ \widehat{\gamma} : s\in \mathbb{R} \mapsto (s,\gamma(s))\in \mathbb{R} \times N $ 
is a solution of the system of differential equations
\begin{equation*}
\left\{
    \begin{aligned}
    \frac{dx_i}{ds} &= \eta_i(t, x), \\
    \frac{dt}{ds}   &= 1,
    \end{aligned}
\right. \qquad i = 1, \dots, n.
\end{equation*} 
The reparametrization $t = t(s) $ shows that $\gamma(t)$ is a solution to
\begin{equation} \label{nonauto}
\frac{dx_i}{dt} = \eta_i(t, x), \qquad i = 1, \dots, n.
\end{equation}
System \eqref{nonauto} is the {\it associated system} with $ X $. Also, a first-order system of ODEs in normal form \eqref{nonauto} gives rise to a $t$-dependent vector field on $N$ of the form 
\[ X(t, x) = \sum_{i = 1}^n \eta_i(t, x) \frac{\partial}{\partial x_i} \]
whose integral curves are of the form $t\mapsto (t,\gamma(t))$, where $\gamma(t)$ is a particular solution to \eqref{nonauto}. This fact justifies identifying $X$ with the t-dependent first-order system of ordinary differential equation (\ref{nonauto}).


For our purposes it is important to relate $t$-dependent vector fields to Lie algebras. A {\it Lie algebra} is a pair $ (V, [\cdot \, , \cdot]) $, where $ V $ is a vector space and $ [\cdot \, , \cdot] : V \times V \to V $ is a bilinear and antisymmetric map that satisfies the Jacobi identity.
The {\it minimal Lie algebra}, $ {\rm Lie}(\mathcal{B}, V, [\cdot \, , \cdot]) $, of a subset $\mathcal{B}\subset V$ of a Lie algebra $(V,[\cdot\,,\cdot])$ is the smallest Lie subalgebra (in the sense of inclusion) in $ V $ that contains  $ \mathcal{B}$. If it does not lead to misunderstanding, $ {\rm Lie}(\mathcal{B}, V, [\cdot \, , \cdot]) $ will simply be denoted by $ {\rm Lie}(\mathcal{B}) $. 
Given a $t$-dependent vector field $ X $ on $N$, we call {\it minimal Lie algebra} of $ X $ the smallest Lie algebra, $V^X$, of vector fields on $N$  that contains all the vector fields $ \{X_t\}_{t \in \mathbb{R}} $. 

\subsection{Lie groups and matrix Lie groups} \label{fundaLie}

Let $G$ be a Lie group and let $e$ be its neutral element. Every  $ g \in G $ defines a right-translation $R_g : h\in G \mapsto hg\in G $ and a left-translation $L_g:h\in G\mapsto gh\in G$ on $G$.
A vector field, $ X^\mathrm{R} $, on $ G $ is {\it right-invariant} if $ X^\mathrm{R}(hg) = R_{g*,h} X^\mathrm{R}(h) $ for every $ h, g \in G $, where $ R_{g*, h} $ is tangent map to  $ R_g $ at $h\in G$. The value of a right-invariant vector field, $ X^\mathrm{R} $, at every point of $ G $ is determined by its value at $e$, since, by definition, $ X^\mathrm{R}(g) = R_{g*,e} X^\mathrm{R}(e) $ for every $ g \in G$. Hence, each right-invariant vector field $X^R$ on $ G $ gives rise to a unique $X^R(e)\in T_eG$ and vice versa. Then, the space of right-invariant vector fields on $G$ is a finite-dimensional Lie algebra.
Similarly, one may define {left-invariant} vector fields on $G$, establish a Lie algebra structure on the space of left-invariant vector fields and set an isomorphism between the space $\mathfrak{g}$ of left-invariant vector fields on $G$ and $T_eG$. The Lie algebra of left-invariant vector fields on $G$, with Lie bracket $ [\cdot \, , \cdot] : \mathfrak{g} \times \mathfrak{g} \to \mathfrak{g} $, induces in $T_eG$ a Lie algebra via the identification of left-invariant vector fields and their values at $e$. Note that we will frequently identify $\mathfrak{g}$ with $T_eG$ to simplify the terminology.

There is a natural mapping from $\mathfrak{g}$ to $G$, the so-called {\it exponential map}, of the form $ \exp :a\in  \mathfrak{g} \mapsto \gamma_a(1)\in G $, where $ \gamma_a: \mathbb{R} \to G $ is the integral curve of the right-invariant vector field $ X^\mathrm{R}_a $ on $G$ satisfying $ X^\mathrm{R}_a(e) = a $ and $ \gamma(0) = e $. If $\mathfrak{g}=\mathfrak{gl}(n,\mathbb{K})$, where $\mathfrak{gl}(n,\mathbb{K})$ is the Lie algebra of $n\times n$ square matrices with entries in a field $\mathbb{K}$ relative to the Lie bracket given by the commutator of matrices, then $\mathfrak{gl}(n,\mathbb{K})$ can be considered as the Lie algebra of the Lie group ${\rm GL}(n,\mathbb{K})$ of $n\times n$ invertible matrices with entries in $\mathbb{K}$. It can be proved that in this case $\exp:X\in \mathfrak{gl}(n,\mathbb{K})\mapsto \exp(X)\in {\rm GL}(n,\mathbb{K})$ retrieves the standard expression of the exponential of a matrix \cite{lee00}, namely
\[ \exp(X) = {\rm I}_n + X + \frac{X^2}{2} + \frac{X^3}{6} + \cdots = \sum_{k = 0}^{\infty} \frac{X^k}{k!},\]
where ${\rm I}_n$ stands for the $n\times n$ identity matrix.

From the definition of the exponential map $\exp:T_eG\rightarrow G$, it follows that $\exp(sa)=\gamma_a(s)$ for each $ s \in \mathbb{R} $ and $a\in T_eG$. Let us show this. Indeed, given the right-invariant vector field $X^\mathrm{R}_{sa} $, where $ sa \in T_eG $,  then 
$$
X^\mathrm{R}_{sa}(g) = R_{g*,e} X^\mathrm{R}_{sa}(e) = R_{g*,e}(sa) = s R_{g*,e}(a),\qquad \forall g\in G.
$$
In particular for $s=1$, it follows that $ X^\mathrm{R}_a(g) =R_{g*,e}(a) $ and, for general $s$, it follows that $X^\mathrm{R}_{sa} = s X^\mathrm{R}_a $. Hence, if $ \gamma_a,\gamma_{sa}:\mathbb{R}\rightarrow G$ are the integral curves of $ X^\mathrm{R}_a $ and $X^\mathrm{R}_{sa}$ with initial condition $ e $, respectively, then it can be proved that, for $u=ts$, one has that
$$
\frac{d}{dt}\gamma_a(ts)=s\frac{d}{du}\gamma_a(u)=sX^{\rm R}_a(\gamma_a(ts)).
$$
and $t\mapsto \gamma_a(st)$ is the integral curve of $X_{sa}^{\rm R}$ with initial condition $e$. Hence, $\gamma_{a}(st)=\gamma_{sa}(t)$.
Therefore, $ \exp(sa) = \gamma_{sa}(1) = \gamma_a(s)$. It is worth stressing that 
Ado's theorem \cite{ado47} shows that every Lie group admits a matrix representation close to its neutral element. 

The exponential map establishes a diffeomorphism from an open neighborhood $U_\mathfrak{g}$ of $0$ in $T_eG$  and $\exp(U_\mathfrak{g})$. More in detail, every basis $ \mathcal{V} = \{ v_1, \dots,v_r \}$ of $T_eG$ gives rise to the so-called  {\it canonical coordinates of the second-kind} related to $\mathcal{V}$ defined by the local diffeomorphism
\begin{equation*}
\begin{array}{ccc}
U_\mathfrak{g}\subset T_eG                & \longrightarrow    & \exp(U_\mathfrak{g}) \subset G \\
(\lambda_1, \dots, \lambda_r) & \mapsto & \prod_{\alpha = 1}^{r} \exp( \lambda_\alpha v_\alpha ) \, ,
\end{array}
\end{equation*}
for an appropriate open neighborhood $U_\mathfrak{g}$  of $0$ in $T_eG\simeq\mathfrak{g}$. 

In matrix Lie groups right-invariant vector fields take a simple useful form . In fact, let $G$ be a matrix Lie group. It can be then considered as a Lie subgroup of ${\rm GL}(n,\mathbb{K})$. Moreover, it can be proved that $T_AG$, for any $A\in G$, can be identified with the space of $n\times n$ square matrices $\mathcal{M}_n(\mathbb{K})$.

Since $ R_A : B\in G \mapsto BA\in  G $, then $ R_{A*,e} (M) = M A \in T_A G,$ for all $M \in T_eG , $ and $A\in {\rm GL}(n,\mathbb{K})$. As a consequence, if $X^{\rm R}(e)=M$ at the neutral element $e$, namely the identity ${\rm I}$, of the matrix Lie group $G$, then
$X^\mathrm{R}(A) = R_{A*,{\rm I}}(X^\mathrm{R}({\rm I})) = R_{A*,{\rm I}}(M) = MA$. It follows that, at any $ A \in G $, every tangent vector $ B \in T_A G $ can be written as $ B = CA $ for a unique $ C \in T_IG$  \cite{Hall,MLC}. 

Let us describe some basic facts on Lie group actions on manifolds induced by Lie algebras of vector fields. It is known that every finite-dimensional Lie algebra, $V$, of vector fields on a manifold $N$ gives rise to a (local) Lie group action
\begin{equation}\label{Accion}
\varphi : G \times N \to N,
\end{equation}
whose fundamental vector fields are given by the elements of $V$ and $G$ is a connected and simply connected Lie group whose Lie algebra is isomorphic to $V$. If the vector fields of $V$ are complete, then the Lie group action (\ref{Accion}) is globally defined on $G\times N$. Let us show how to obtain $\varphi$ from $V$, which will be of crucial importance in this work.

Let us restrict ourselves to an open neighborhood $ U_G $ of the neutral element of $G$, where we can use canonical coordinates of the second-kind related to a basis $ \{ v_1, \dots,v_r \} $ of $\mathfrak{g}$. Then, each $ g \in U_G $ can be expressed as
\begin{equation} \label{ccse}
g = \prod_{\alpha = 1}^{r} \exp( \lambda_\alpha v_\alpha ) , 
\end{equation}
for certain uniquely defined parameters $ \lambda_1, \dots, \lambda_r \in \mathbb{R} $. To determine $\varphi$, we determine the curves
\begin{equation} \label{curv_accion}
\gamma_x^\alpha : \mathbb{R} \to N : t \mapsto \varphi(\exp(t v_\alpha), x), \qquad \alpha = 1, \dots, r,
\end{equation}
where $\gamma^\alpha_x$ must be the integral curve of $X_\alpha$ for $\alpha=1,\ldots,r$. Indeed, for any element $ g \in U_G \subset G $ expressed as in \eqref{ccse}, using the intrinsic properties of a Lie group action, 
\begin{equation*}
\varphi(g, x) = \varphi\left( \prod_{\alpha = 1}^{r} \exp( \lambda_\alpha v_\alpha ), x \right)=\left( \varphi( \exp( \lambda_1 v_1 ))\cdot \varphi ( \exp(\lambda_2 v_2))\cdot \varphi(\exp(\lambda_r v_r)),x \right) ,
\end{equation*} 
the action is completely defined for any $ g \in U_G \subset G $.

In this work we will deal with some particular matrix Lie groups, starting from
the general linear matrix group $ \mathrm{GL}(n, \mathbb{K}) $, where we recall that $\mathbb{K}$ may be $\mathbb{R}$ or $\mathbb{C}$. 
As it is well known, any closed subgroup of $\mathrm{GL}(n, \mathbb{K})$ is also a matrix Lie group \cite[Theorem 15.29, pg. 392]{lee00}. In the forthcoming pages we will work with some of those subgroups such as $ \mathrm{SL}(n, \mathbb{R})$, the Lie group formed by $n\times n$ real matrices with unit determinant. Moreover, for future reference we recall that the Lie algebra of ${\rm SL}(n,\mathbb{R})$, i.e., $\mathfrak{sl}(n, \mathbb{R})$, is the space of traceless $n\times n$ real matrices  \cite{Hall, Sattinger}.

\subsection{Lie systems} \label{capLie}

The Lie Theorem \cite{carinena07} states that a Lie system is a $t$-dependent system of (first-order) ordinary differential equations that describes the integral curves of a $t$-dependent vector field that takes values in a finite-dimensional Lie algebra of vector fields, namely the aforementioned Vessiot-Guldberg Lie algebra (VG) \cite{carinena00, lucas11}. As we also mentioned previously, one of the most important characteristics of Lie systems is that they admit (generally nonlinear) superposition rules and a plethora of mathematical properties mediated by the {Lie theorem} \cite{carinena07}. Furthermore, some Lie systems can be studied via a Hamiltonian formulation \cite{carinena11, araujo20}.

In this section we introduce some of these fundamental concepts in the theory of Lie systems. In this way, we start by introducing solutions of Lie systems in terms of superposition rules.

On a first approximation, a {\it Lie system} is a first-order system of ODEs that admits a superposition rule.

\begin{definition} 
A {\it superposition rule} for a system $ X $ on $ N $ is a map $ \Phi : N^{m} \times N \to N $ such that the general solution $ x(t) $ of $ X $ can be written as $ x(t) = \Phi(x_{(1)}(t), \dots, x_{(m)}(t); \rho) $, where $ x_{(1)}(t), \dots, $ $ x_{(m)}(t) $ is a generic family of particular solutions and $ \rho $ is a point in $ N $ related to the initial conditions of $X$.
\end{definition}

A classic example of Lie system is the Riccati equation  \cite[Example 3.3]{araujo20}, that is,
\begin{equation} \label{Ric}
\frac{dx}{dt} = b_1(t) + b_2(t) x + b_3(t) x^2,\qquad x\in \mathbb{R},
\end{equation} 
with $ b_1(t), b_2(t), b_3(t) $ being arbitrary functions of $t$. It is known then that the general solution, $x(t)$, of the Riccati equation can be written as
\begin{equation} \label{SupRiccati}
x(t) = \frac{x_{(2)}(t)(x_{(3)}(t) - x_{(1)}(t)) + \rho x_{(3)}(t) (x_{(1)}(t) - x_{(2)}(t))}{(x_{(3)}(t) - x_{(1)}(t)) + \rho (x_{(1)}(t) - x_{(2)}(t))},
\end{equation}
where $ x_{(1)}(t), x_{(2)}(t), x_{(3)}(t) $ are three different particular solutions of \eqref{Ric} and $ \rho \in \mathbb{R} $ is an arbitrary constant. This implies that the Riccati equation admits a superposition rule $\Phi:\mathbb{R}^3\times \mathbb{R}\rightarrow \mathbb{R}$ such that
$$
\Phi(x_{(1)},x_{(2)},x_{(3)},\rho)=\frac{x_{(2)}(x_{(3)} - x_{(1)}) + \rho x_{(3)} (x_{(1)} - x_{(2)})}{(x_{(3)} - x_{(1)}) + \rho (x_{(1)} - x_{(2)})}.
$$

The conditions that guarantee the existence of a superposition rule are gathered in the  Lie theorem \cite[Theorem 44]{lie93}.  

\begin{theorem}[Lie theorem] A first-order system $ X $ on $ N $,
\begin{equation} \label{genLiesys1}
\frac{dx}{dt} = X(t, x), \qquad  x \in N, \qquad  X\in\mathfrak{X}_t(N),
\end{equation}
admits a superposition rule if and only if $ X $ can be written as
\begin{equation} \label{genLiesys}
X(t,x) = \sum_{\alpha = 1}^r b_\alpha(t) X_\alpha(x),\qquad t\in \mathbb{R},\qquad x\in N,
\end{equation}
for a certain family $ b_1(t), \dots, b_r(t) $ of $t$-dependent functions and a family of vector fields $ X_1, \dots, $ $ X_r $ on $N$ that generate an $r$-dimensional Lie algebra of vector fields.
\end{theorem}

The Lie theorem yields  that 
every Lie system $ X $ is related to (at least) one VG Lie algebra, $V$, that satisfies that Lie($\{ X_t \}_{t \in \mathbb{R}}$) $\subset V $. This implies that the minimal Lie algebra has to be finite-dimensional, and vice versa \cite{lucas11}. 


\begin{example}
The  $t$-dependent vector field on the real line associated with \eqref{Ric} is $ X = b_1(t) X_1 + b_2(t) X_2 + b_3(t) X_3 $, where $X_1,X_2,X_3$ are vector fields on $\mathbb{R}$ given by
\begin{equation*}
X_1 = \frac{\partial}{\partial x}, \qquad X_2 = x \frac{\partial}{\partial x}, \qquad X_3 = x^2 \frac{\partial}{\partial x}.
\end{equation*}
Since the commutation relations are
\begin{equation}
    [X_1,X_2]=X_1,\quad [X_1,X_3]=2X_2,\quad [X_2,X_3]=X_3,
\end{equation}
the vector fields $ X_1, X_2, X_3 $ generate a VG Lie algebra isomorphic to $ \mathfrak{sl} (2, \mathbb{R}) $.
Then, the Lie theorem guarantees that \eqref{Ric} admits a superposition rule, which is precisely the one shown in \eqref{SupRiccati}.
\end{example}

\subsubsection{Automorphic Lie systems} \label{pasoNtoG}

The general solution of a Lie system on $N$ with a VG Lie algebra, $V$, can be obtained from a single particular solution of a Lie system on a Lie group $G$ whose Lie algebra is isomorphic to $V$, a so-called {\it automophic Lie system} \cite[\S 1.4]{lucas11}. As the automorphic Lie system notion is going to be central in our paper, let us study  it in some detail (see \cite{lucas11} for details).

\begin{definition} An {\it  automorphic Lie system} is a $t$-dependent system of first-order differential equations on a Lie group $G$ of the form
\begin{equation}\label{rhshat}
\frac{dg}{dt}=\sum_{\alpha=1}^rb_\alpha(t)X_\alpha^R(g),\qquad g\in G,\quad t\in\mathbb{R},
\end{equation}
where $\{X_1^R,\ldots,X_r^R\}$ is a basis of the space of right-invariant vector fields on $G$ and $b_1(t),\ldots,b_r(t)$ are arbitrary $t$-dependent functions.  Furthermore, we shall refer to the right-hand side of equation \eqref{rhshat} as $\widehat{X}^G_R(t,g)$, i.e., $\widehat{X}^G_R(t,g)=\sum_{\alpha = 1}^{r} b_\alpha(t) X_\alpha^\mathrm{R}(g)$.
\end{definition}

 Because of right-invariant vector fields, systems in the form of
$ \widehat{X} $ have the following important property.
\begin{proposition}\label{Prop:RA} (See \cite[\S 1.3]{lucas11})
Given a Lie group $ G $ and a particular solution $ g(t) $ of the Lie system defined on $ G $, as 
\begin{equation}\label{LiesysLiegroup}
\frac{dg}{dt} = \sum_{\alpha = 1}^{r} b_\alpha(t) X_\alpha^\mathrm{R}(g)=\widehat{X}^G_R(t,g), \, 
\end{equation}
where $ b_1(t),\ldots,b_r(t) $ are arbitrary $t$-dependent functions and $ X_1^\mathrm{R},\ldots,X_r^\mathrm{R}$ are right-invariant vector fields, we have that $ g(t) h $ is also a solution of \eqref{LiesysLiegroup} for each $ h \in G $. 
\end{proposition}

An immediate consequence of Proposition \ref{Prop:RA} is that, once we know a particular solution of  $ \widehat{X}^G_R $, any other solution can be obtained simply by  multiplying the known solution on the right by any element in $ G $. More concretely, if we know a solution
 $ g(t) $ of (\ref{LiesysLiegroup}), then the solution $ h(t) $ of (\ref{LiesysLiegroup}) with initial condition $ h(0) = g(0)h_0 $ can be expressed as $h(t)=g(t) h_0 $. This justifies that henceforth we only worry about finding one particular solution $g(t)$ of $ \widehat{X}^G_R $, e.g. the one that fulfills $ g(0) = e $. The previous result can be understood in terms of the Lie theorem or via superposition rules. In fact, since \eqref{LiesysLiegroup} admits a superposition rule $\Phi:(g,h)\in G\times G\mapsto gh\in G$, the system (\ref{Prop:RA})   must be a Lie system. Alternatively, the same result follows from the Lie Theorem and the fact that the right-invariant vector fields on $G$ span a finite-dimensional Lie algebra of vector fields.

There several reasons to study automorphic Lie systems. One is that they can be locally written around the neutral element of its Lie group in the form 
$$
\frac{dA}{dt}=B(t)A,\qquad A\in {\rm GL}(n,\mathbb{K}),\quad B(t)\in \mathcal{M}_n(\mathbb{K}),
$$ 
where $\mathcal{M}_n(\mathbb{K})$ is the set of $n\times n$ matrices o coefficients in $\mathbb{K}$, for every $t\in \mathbb{R}$.

The main reason to study automorphic Lie systems is given by the following results, which show how they can be used to solve any Lie system on a manifold. Let us start with a Lie system $ X $ defined on $ N $. Hence, $X$ can be written as
\begin{equation} \label{Lie-Scheffer}
\frac{dx}{dt} = \sum_{\alpha = 1}^{r} b_\alpha(t) X_\alpha ,
\end{equation}
for certain $t$-dependent functions $ b_1(t),\ldots,b_r(t) $ and vector fields $ X_1,\ldots,X_r\in  \mathfrak{X}(N) $ that generate a $r$-dimensional  dimensional VG Lie algebra.
The VG Lie algebra $V$ is always isomorphic to the Lie algebra $ \mathfrak{g} $ of a certain  Lie group $G$. The VG Lie algebra spanned by $X_1,\ldots,X_r$ gives rise to a (local) Lie group action $\varphi:G\times N\rightarrow N$ whose fundamental vector fields are those of $V$. In particular, there exists a basis $\{v_1,\ldots,v_r\}$ in $\mathfrak{g}$ so that
$$
\frac{d}{dt}\bigg|_{t=0}\varphi(\exp(tv_\alpha),x)=X_\alpha(g),\qquad \alpha=1,\ldots,r.
$$
In other words, $\varphi_\alpha:(t,x)\in \mathbb{R}\times N\mapsto \varphi(\exp(tv_\alpha),x)\in N$ is the flow of the vector field $X_\alpha$ for $\alpha=1,\ldots,r$. Note that if $[X_\alpha,X_\beta]=\sum_{\gamma=1}^rc_{\alpha\beta}^\gamma X_\gamma$ for $\alpha,\beta=1,\ldots,r$, then $[v_\alpha,v_\beta]=-\sum_{\gamma=1}^rc_{\alpha\beta}^\gamma v_\gamma$  for $\alpha,\beta=1,\ldots,r$ (cf. \cite{carinena00}).



To determine the exact form of the Lie group action $ \varphi : G \times N \to N $ as in \eqref{curv_accion}, we impose
\begin{equation} \label{def_accion}
\varphi(\exp(\lambda_\alpha v_\alpha), x) = \varphi_\alpha(\lambda_\alpha, x) \qquad \forall \, \alpha = 1, \dots, r,\qquad \forall x\in N,
\end{equation}
where $ \lambda_1,\ldots,\lambda_r \in \mathbb{R} $. While we stay in a neighborhood $U$ of the origin of $G$, where every element $g\in U$ can be written in the form 
$$
g=\exp(\lambda_1v_1)\cdot\ldots\cdot \exp(\lambda_rv_r),
$$
then the relations \eqref{def_accion} and the properties of $\varphi$ allow us to determine $ \varphi $ on $U$. If we fix $ x \in N $, the right-hand side of the equality turns into an integral curve of the vector field $ X_\alpha $, this is why (\ref{def_accion}) holds.

\begin{proposition} \label{prop.accion} (see \cite{carinena00,lucas11} for details) Let $ g(t) $ be a solution to the system 
$$
\frac{dg}{dt}=\sum_{\alpha=1}^rb_\alpha(t)X^R(g),\qquad \forall t\in \mathbb{R},\quad g\in G.
$$
Then, $ x(t) = \varphi(g(t), x_0) $ is a solution of $ X=\sum_{\alpha=1}^rb_\alpha(t)X_\alpha $, where $ x_0 \in N $. In particular, if one takes the solution $ g(t) $ that satisfies the initial condition $ g(0) = e $, then $ x(t) $ is the solution of $ X $ such that $ x(0) = x_0 $. 
\end{proposition}

Let us study a particularly relevant form of automorphic Lie systems that will be used hereafter. If $\mathfrak{g}$ is a finite-dimensional Lie algebra, then Ado's theorem \cite{ado47} guarantees that $\mathfrak{g}$ is isomorphic to a matrix Lie algebra $\mathfrak{g}_M$. Let $ \mathcal{V} = \{M_1, \ldots, M_r\} $ be a basis of $ \mathfrak{g}_M\subset\mathcal{M}_n(\mathbb{R})$. 
As reviewed in Section \ref{fundaLie},  each  $ M_\alpha $ gives rise to a right-invariant vector field $X^R_\alpha(g)=M_\alpha g$, with $g\in G$, on $ G $. These vector fields have the opposite  commutation relations than the (matrix) elements of the basis. 

In the case of matrix Lie groups, the system  \eqref{LiesysLiegroup} takes a simpler form. Let $ Y(t) $ be the matrix associated with the element $ g(t)\in G $. Using the right invariance property of each $ X_\alpha^\mathrm{R} $, we have that
\begin{equation*}
\frac{dY}{dt} = \sum_{\alpha = 1}^{r} b_\alpha(t) X_\alpha^\mathrm{R}(Y(t)) = \sum_{\alpha = 1}^{r} b_\alpha(t) R_{Y(t)*, e} \left( X_\alpha^\mathrm{R}(e) \right) = \sum_{\alpha = 1}^{r} b_\alpha(t) R_{Y(t)*, e} (M_\alpha) .
\end{equation*}
We can write the last term as
\[ \sum_{\alpha = 1}^{r} b_\alpha(t) R_{Y(t)*, e} (M_\alpha) = \sum_{\alpha = 1}^{r} b_\alpha(t) M_\alpha Y(t) , \]
in such a way that for matrix Lie groups, the system on the Lie group is
\begin{equation} \label{sistemaGM}
\frac{dY}{dt} = A(t) Y(t) , \qquad Y(0) = I , \qquad \text{with} \quad  A(t) = \sum_{\alpha = 1}^{r} b_\alpha(t) M_\alpha ,
\end{equation}
where $ I $ is the identity matrix (which corresponds with the neutral element of the matrix Lie group) and the matrices $ M_\alpha $ form a finite-dimensional Lie algebra, which is anti-isomorphic to the VG Lie algebra of the system (by anti-isomorphic we imply that the systems have the same constants of structure but that they differ in one sign).

There exist various methods to solve system \eqref{LiesysLiegroup} analytically \cite[\S 2.2]{sardon15}, such as the Levi decomposition \cite{levi05} or the theory of reduction of Lie systems \cite[Theorem 2]{carinena01}. In some cases, it is relatively easy to solve it,  as is the case where $ b_1,\ldots,b_r $ are constants. We will depict an example in this particular case in Section \ref{capSLN}. 
Nonetheless, we are interested in a numerical approach, since we will try to solve the automorphic Lie system  with adapted geometric integrators. The solutions on the Lie group can be straightforwardly translated into solutions on the manifold for the Lie system defined on $N$ via  the Lie group action \eqref{Accion}.

To finish this section, we will employ the previous developments in order to define our novel procedure to (geometrically) construct a continuous solution of a given Lie system.


\medskip

{\bf The 7 step method}: Reduction procedure to automorphic Lie system
\label{Metodo7Pasos}

The method can be itemized in the following seven steps:

\begin{enumerate}

\item  We identify the VG Lie algebra of vector fields $ X_1, \dots, X_r $  that defines the Lie system on $N$.
\item  We look for a Lie algebra  {\rm $\mathfrak{g}$} isomorphic to the VG Lie algebra, whose basis is $\{M_1, \dots, M_r \}\in \mathcal{M}_{n\times n}(\mathbb{R}) $ with the same structure constants of $X_1,\ldots,X_r$ in absolute value, but with a negative sign.
\item We integrate the vector fields $ X_1,\ldots,X_r $ to obtain their respective flows $ \Phi_\alpha : \mathbb{R} \times N \to N $ with $\alpha=1,\ldots,r$.
\item Using canonical coordinates of the second kind and the previous flows we construct the Lie group action  $ \varphi : G \times N \to N $ using expressions \eqref{def_accion}. 
\item We define an automorphic Lie system $ \widehat{X}^G_R $ on the Lie group $ G $ associated with $ \mathfrak{g} $ as in \eqref{LiesysLiegroup}.
\item We compute the solution of the system $\widehat{X}^G_R$ that fulfils $ g(0) = e $. 
\item Finally, we retrieve the solution for $ X $ on $N$ through the expression $ x(t) = \varphi(g(t), x_0) $. 
\end{enumerate}

\section{Discretization of Lie systems} \label{cap.discre}

This section adapts known numerical methods on Lie groups to automorphic Lie systems. For this purpose, we start by reviewing briefly some fundamentals on numerical methods for ordinary differential equations and Lie groups \cite{hairer93, isaacson66, quarteroni07}, and later focus on two specific numerical methods on Lie groups, the Magnus expansion and RKMK methods \cite{iserles05,iserles99,  munthe98, munthe99, Zanna}. 

Recall that, in this paper, we focus on ordinary differential equations of the form
\begin{equation}\label{DynSysNumer}
\frac{dx}{dt} = f(t, x), \qquad t\in[a,b], \qquad  x(t) \in N, \qquad  f\in\mathfrak{X}_t(N).
\end{equation}
When $N$ is (or diffeomorphic to) an  Euclidean space, there is a plethora of numerical schemes approximating the analytic solution $x(t)$ of (\ref{DynSysNumer}) \cite{hairer93, isaacson66}.  We will focus on {\it one-step methods with fixed time step}. By that we mean that solutions  are approximated by a sequence of numbers $x_k=x(t_k)\in N$ with $t_k=a+kh$, $h=(b-a)/\mathcal{N}$, $b>a$ and
\begin{equation}\label{OneStepMethod}
\frac{x_{k+1}-x_k}{h}=f_h(t_k,x_k,x_{k+1}),
\end{equation}
where $\mathcal{N}$ is the number of steps our time interval is divided to.  We call $h$ the {\it time step}, which is fixed, while  $f_h:\mathbb{R}\times N\times N\rightarrow TN$ is a discrete vector field, which (recall that, for now, we set $N$ to be a Euclidean space with norm $\|\cdot\|$) is a given approximation of $f$
in \eqref{DynSysNumer}. As usual, we shall denote the local truncation error by $E_h$, where
\begin{equation}\label{LTE}
E_h=||x_{k+1}-x(t_{k+1})||, 
\end{equation}
and say that the method is of order $r$ if $E_h=\mathcal{O}(h^{r+1})$ for $h\rightarrow 0$, i.e. $\lim_{h\rightarrow0}|E_h/h^{r+1}|<\infty$. Regarding the global error
\[
E_\mathcal{N}=||x_\mathcal{N}-x(b)||,
\]
we shall say that the method is \textit{convergent} of order $r$ if $E_\mathcal{N}=\mathcal{O}(h^{r})$, when $h\rightarrow 0$. As for the simulations, we pick the following norm in order to define the global error, that is
\[
E_\mathcal{N}=\max_{k=1,\ldots,\mathcal{N}}||x(t_k)-x_k||.
\]

Given the relevant examples in this paper, e.g., Ricatti equations, where $N=\mathbb{R}^n$, we will employ classical methods to approximate \eqref{DynSysNumer}, particularly the Heun method (convergent of order 2) and RK4 (convergent of order 4), and compare to our novel discretization proposal.
\color{black}

\subsection{Numerical methods on matrix Lie groups}
\label{LieGroupsMethods}

Our purpose is to numerically solve the initial condition problem for system \eqref{sistemaGM} defined on a matrix Lie group $ G $ of the form
\begin{equation} \label{eq.grupo}
\frac{dY}{dt} = A(t) Y \qquad \text{with} \qquad Y(0) = I, 
\end{equation}
where $ Y \in G $ while $ A(t) \in \mathfrak{g} \cong T_e G $ is a given $t$-dependent matrix and $ I $ is the identity matrix in $G$. That is, we are searching for a discrete sequence $\{Y_k\}_{k=0,\ldots,\mathcal{N}}$ such that $Y_k\in G$. In a neighborhood of the zero in $T_eG$, the exponential map defines a diffeomorphism onto an open subset of the neutral element of $G$ and the problem is equivalent to searching for a curve
 $ \Omega(t)$ in $\mathfrak{g} $ such that
\begin{equation} \label{eq.omega}
Y(t) = \exp(\Omega(t)) . 
\end{equation}
This ansatz helps us to transform \eqref{eq.grupo}, which is defined in a nonlinear space, into a new problem in a linear space, namely the Lie algebra $\mathfrak{g}\simeq T_eG$. This is expressed in the classical result by Magnus \cite{Magnus}.

\begin{theorem}[Magnus, 1954]

The solution of the matrix Lie group \eqref{eq.grupo} in $G$ can be written for values of $t$ close enough to zero, as $Y(t)=\exp(\Omega(t))$, where $\Omega(t)$ is the solution of the initial value problem
\begin{equation} \label{eq.algebra} 
\frac{d \Omega}{dt} = \operatorname{dexp}^{-1}_{\Omega(t)}(A(t)),\qquad\quad \Omega(0) = {\bf 0} \, ,
\end{equation}
where $ {\bf 0} $ is the zero element in $T_eG$.
\end{theorem}

When we are dealing with matrix Lie groups and Lie algebras, the $\mbox{dexp}^{-1}$ is given by
\begin{equation} \label{eq.algebraS}
 \mbox{dexp}^{-1}_{\Omega}(H)= \sum_{j = 0}^{\infty} \frac{B_j}{j!} \operatorname{ad}_{\Omega}^j(H), \, 
\end{equation} 
where the $\{B_j\}_{j=0,\ldots,\infty}$ are the Bernoulli numbers and $\mbox{ad}_{\Omega}(H)=[\Omega,H]=\Omega\,H-H\,\Omega.$ The convergence of the series \eqref{eq.algebraS} is ensured as long as a certain convergence condition is satisfied \cite{Magnus}. 

If we try to integrate \eqref{eq.algebra} applying a numerical method directly (note that, now, we could employ one-step methods \eqref{OneStepMethod} safely), $ \Omega(t) $  might sometimes drift too much away from the origin and the exponential map would not work. This would be a problem, since we are assuming that $ \Omega(t) $ stays in a neighborhood of the origin of $\mathfrak{g}$ where the exponential map defines a local diffeomorphism with the Lie group. Since we still do not know how to characterize this neighborhood, it is necessary to adopt  a strategy  that allows us to resolve \eqref{eq.algebra} sufficiently close to the origin. The thing to do is to change the coordinate system in each iteration of the numerical method. In the next lines we explain how this is achieved. 

Consider now the restriction of the exponential map given by
\begin{align*} 
\exp : U_{\mathfrak{g}} \subset \mathfrak{g} &\to \exp(U_\mathfrak{g}) \subset G ,\\
A &\mapsto \exp(A) 
\end{align*}
so that this map establishes a diffeomorphism between
an open neighborhood $ U_{\mathfrak{g}} $ around the origin in $ \mathfrak{g} $ and its image. Since the elements of the matrix Lie group are invertible matrices, the map $ U_{\mathfrak{g}} \to \exp(U_\mathfrak{g})Y_0 \subset G : A \mapsto \exp(A) Y_0 $ from $ U_{\mathfrak{g}} \subset \mathfrak{g} $ to the set
\[ \exp(A)Y_0 = \{ Y \in G \ : \exists X\in U_\mathfrak{g}, Y = X Y_0\} \] is also a diffeomorphism. This map gives rise to the so-called {\it first-order canonical coordinates centered} at $ Y_0 $. 

As well-known, the solutions of \eqref{eq.algebra}  are curves in $\mathfrak{g}$ whose images by the exponential map are solutions to \eqref{eq.grupo}. In particular, the solution $ \Omega^{(0)}(t) $ of system \eqref{eq.grupo} such that $\Omega^{(0)}(0)$ is the zero matrix in $T_{\rm Id}G$, namely ${\bf 0}$,  
corresponds with the solution $Y^{(e)}(t)$ of the system on $G$ such that
$Y^{(e)}(0)=I$. Now, for a certain $ t = t_k $, the solution $ \Omega^{(t_k)}(t)$ in $\mathfrak{g} $ such that $ \Omega^{(t_k)}(t_k) = {\bf 0}$, corresponds with $ Y^{(e)}(t) $
via first-order canonical coordinates centered at
 $ Y^{(e)}(t_k) \in G $, since 
\[ \exp(\Omega^{(t_k)}(t_k)) Y^{(e)}(t_k) = \exp({\bf 0}) Y^{(e)}(t_k) = Y^{(e)}(t_k) , \]
and the existence and uniqueness theorem guarantees $ \exp(\Omega^{(0)}(t)) = \exp(\Omega^{(t_k)}(t)) Y^{(e)}(t_k) $ around $ t_k $. In this way, we can use the curve
 $ \Omega^{(t_k)}(t)$ and the canonical coordinates centered on  $ Y^{(e)}(t_k) $ to obtain values for the solution of \eqref{eq.grupo} in the proximity of $ t = t_k $, instead of using $ \Omega^{(0)}(t) $. Whilst the curve $ \Omega^{(0)}(t) $ could be far from the origin of coordinates for $t_k$,  we know that $ \Omega^{(t_k)}(t) $ will be close, by definition. Applying this idea in each iteration of the numerical method, we are changing the curve in $ \mathfrak{g} $ to obtain the approximate solution of \eqref{eq.grupo} while we stay near the origin (as long as the time step is small enough). 
 
 Thus, what is left is defining proper numerical methods for \eqref{eq.algebra} whose solution, i.e. $\{\Omega_k\}_{k=0,\ldots,\mathcal{N}}$, via the exponential map, provides us with a numerical solution of \eqref{eq.grupo} remaining in $G$. In other words, the general Lie group method defined this way \cite{iserles99, iserles05} can be set by the recursion
 \begin{equation}\label{GeneralLieGroupMethod}
     Y_{k+1}=e^{\Omega_k}\,Y_k.
 \end{equation}
Next, we introduce two relevant families of numerical methods providing $\{\Omega_k\}_{k=0,\ldots,\mathcal{N}}$.

\subsubsection{The Magnus method}
Based on the work by Magnus, the Magnus method was introduced in \cite{iserles99,iserles98}. The starting point of this method is to resolve equation \eqref{eq.algebra} by means of the Picard procedure. This method assures that a given sequence of functions converges to the solution of \eqref{eq.algebra} in a small enough neighborhood. Operating, one obtains the \textit{Magnus expansion}
\begin{equation} \label{ex.magnus}
\Omega(t) = \sum_{k = 0}^{\infty} H_k(t) , 
\end{equation}
where each $ H_k(t)$ is a linear combination of iterated commutators. The first three terms are given by
\begin{align*}
H_0(t) &= \int_0^t A(\xi_1) d \xi_1 \, , \\
H_1(t) &= -\frac{1}{2} \int_0^t\left[ \int_0^{\xi_1} A(\xi_2) d\xi_2, A(\xi_1) \right] d\xi_1 \, , \\
H_2(t) &= \frac{1}{12} \int_0^t\left[ \int_0^{\xi_1} A(\xi_2) d\xi_2, \left[ \int_0^{\xi_1} A(\xi_2) d\xi_2, A(\xi_1) \right] \right] d\xi_1 \\
&\qquad \qquad + \frac{1}{4} \int_0^t\left[ \int_0^{\xi_1} \left[ \int_0^{\xi_1} A(\xi_2) d\xi_2, A(\xi_1) \right] d\xi_2, A(\xi_1) \right] d\xi_1 \, .
\end{align*}
Note that the Magnus expansion
\eqref{ex.magnus} converges absolutely in a given norm for every $ t \geq 0 $ such that \cite[p. 48]{iserles05}
\[ \int_0^t \| A(\xi) \| d\xi \leq \int_{0}^{2 \pi} \frac{d \xi}{4 + \xi [1 - \cot(\xi / 2)]} \approx 1,086868702. \]

In practice, if we work with the Magnus expansion we need a way to handle the infinite series and calculate the iterated integrals. Iserles and N\o rsett proposed a method based on binary trees \cite{iserles99, iserles98}. In \cite[\S 4.3]{iserles05} we can find a method to truncate the series in such a way that one obtains the desired order of convergence. Similarly, \cite[\S 5]{iserles05} discusses in detail how the iterated integrals can be integrated numerically.  
In our case, for practical reasons we will implement the Magnus method following the guidelines of Blanes, Casas \& Ros \cite{blanes00}, which is based on a Taylor series of $ A(t) $ in \eqref{eq.grupo} around the point $ t = h/2 $ (recall that, in the Lie group and Lie algebra equations we are setting the initial time $t_0=a=0$). With this technique one is able to achieve different orders of convergence. In particular, we will use the second and fourth order convergence methods \cite[\S 3.2]{blanes00}, although one can build up to eighth order methods.

The second-order approximation is
\[ \exp(\Omega(h)) = \exp(h a_0) + \mathcal{O}(h^3) \]
and the forth-order one reads
\[ \exp(\Omega(h)) = \exp\left( h a_0 + \frac{1}{12} h^3 a_2 - \frac{1}{12} h^3 [ a_0, a_1 ] \right)  + \mathcal{O}(h^5) , \]
where $ \Omega(0) = {\bf 0} $ and 
\[ a_i = \frac{1}{i!} \left. \frac{d^i}{dt^i} A(t) \right|_{t = h/2} \qquad i = 0, 1, 2 . \]

As we see from the definition, the first method computes the first and second derivative of matrix $ A(t) $. Applying the coordinate change in each iteration \eqref{GeneralLieGroupMethod}, we can implement it through the following equations:
\begin{align} \label{met.mag2}
Y_{k + 1} = \exp \left[ h A \left( t_k + \frac{h}{2} \right) \right]  Y_k.  \qquad \text{[Order 2]} \\[10pt]
\left. \begin{gathered} \label{met.mag4}
           Y_{k + 1} = \exp\left( h a_0 +  h^3 (a_2 - [ a_0, a_1 ]) \right) Y_k , \\[3pt]
           t_{1/2} = t_k + \frac{h}{2} , \quad a_0 = A(t_{1/2} ) , \quad a_1 = \frac{\dot{A} (t_{1/2})}{12} , \quad a_2 = \frac{\ddot{A} (t_{1/2})}{24},       \end{gathered} \ \right\rbrace \qquad \text{[Order 4]}
\end{align}
where $\dot A(t_0),\ddot A(t_0)$ stand for the first and second derivatives of $A(t)$ in terms of $t$ at $t_0$.  
Note that the convergence order is defined for the Lie group dynamics \eqref{eq.grupo}. That is, when we say that the above methods are convergent of order 2, for instance, that means $E_\mathcal{N}=||Y_\mathcal{N}-Y(b)||=\mathcal{O}(h^2)$, with $h\rightarrow 0$, for a proper Lie matrix norm.

\subsubsection{The Runge-Kutta-Munthe-Kaas method}

Changing the coordinate system in each step, as explained in previous sections, the classical RK methods applied to Lie groups give rise to the so-called Runge-Kutta-Munthe-Kaas (RKMK) methods \cite{munthe98, munthe99}. The equations that implement the method are
\begin{align*}
\begin{aligned}
\Theta_j  &= h \sum_{l = 1}^{s} a_{jl} F_l , \\
F_j       &= \operatorname{dexp}_{\Theta_j}^{-1}(A(t_k + c_j h)), \\
\Theta    &= h \sum_{l = 1}^{s} b_l F_l , \\
Y_{k + 1} &= \exp(\Theta) Y_k .
\end{aligned} 
\begin{aligned} \left. 
\vphantom{ \begin{aligned}
              \Theta_j  &= h \sum_{l = 1}^{s} a_{jl} F_l , \\
              F_j       &= \operatorname{dexp}_{\Theta_j}^{-1}(A(t_k + c_j h)),
           \end{aligned}
} \quad \right\rbrace \qquad j = 1, \dots, s , \\
\vphantom{ \begin{aligned}
             \Theta    = h \sum_{l = 1}^{s} b_l F_l , \\
             Y_{k + 1} = \exp(\Theta) Y_k .
           \end{aligned}
}
\end{aligned}  
\end{align*}
where the constants $ \{ a_{jl} \}_{j,l = 1}^s $, $ \{ b_l \}_{l = 1}^s $, $ \{ c_j \}_{j = 1}^s $ can be obtained from a Butcher's table \cite[\S 11.8]{quarteroni07} (note that $s$ is the number of stages of the usual RK methods). Apart from this, we have the consistency condition $ \sum_{l = 1}^s b_l = 1 $.
As the equation that we want to solve comes in the shape of an infinite series, it is necessary to study how we evaluate the function $ \operatorname{dexp}_{\Omega(t)}^{-1} $. For this, we need to use truncated series up to a certain order in such a way that the order of convergence of the underlying classical RK is preserved. If the classical RK is of order $p$ and the truncated series of \eqref{eq.algebra} is up to order $j$, such that
$j \geq p - 2 $, then the RKMK method is of order $p$ (see \cite{munthe98,munthe99} and \cite[Theorem 8.5, p. 124]{hairer06}). Again, this convergence order refers to the equation in the Lie group \eqref{eq.grupo}.

Let us now determine the RKMK method associated with the explicit Runge--Kutta whose Butcher's table is
\begin{center}
\begin{tabular}{c|cccc}
$ 0 $   &         &         &         &         \\
$ 1/2 $ & $ 1/2 $ &         &         &         \\
$ 1/2 $ & $ 0 $   & $ 1/2 $ &         &         \\
$ 1 $   & $ 0 $   & $ 0 $   & $ 1 $   &         \\ \hline
        & $ 1/6 $ & $ 1/3 $ & $ 1/3 $ & $ 1/6 $ 
\end{tabular} 
\end{center}
that is a Runge-Kutta of order 4 (RK4). This implies that we need to truncate the series $\operatorname{dexp}_{\Omega(t)}^{-1} $ at $ j = 2 $:
\begin{equation} \label{dexp.2} 
\operatorname{dexp}_{\Omega}^{-1}(A) \approx A - \frac{1}{2} [\Omega, A] + \frac{1}{12} [\Omega, [\Omega, A]] . 
\end{equation}
Then, the RKMK implementation for the given Butcher's table is
\begin{equation} \label{RKMK4}
\left. \begin{aligned}
F_1  &= \operatorname{dexp}_{O_n}^{-1}(A(t_k)) , \\
F_2  &= \operatorname{dexp}_{\frac{1}{2} h F_1}^{-1} \left( A \left( t_k + \frac{1}{2} h \right) \right) , \\
F_3  &= \operatorname{dexp}_{\frac{1}{2} h F_2}^{-1} \left( A \left( t_k + \frac{1}{2} h \right) \right) , \\
F_4  &= \operatorname{dexp}_{h F_3}^{-1}(A(t_k + h))  ,
\end{aligned} \quad \right\rbrace \quad 
\begin{gathered}
\Theta    = \frac{h}{6} (F_1 + 2 F_2 + 2 F_3 + F_4) , \\
Y_{k + 1} = \exp(\Theta) Y_k ,
\end{gathered}
\end{equation}
where $ \operatorname{dexp}^{-1} $ is \eqref{dexp.2}.

It is interesting to note that the method obtained in the previous section using the Magnus expansion \eqref{met.mag2}
can be retrieved by a RKMK method associated with the following Butcher's table:
\begin{center}
\begin{tabular}{c|cc}
$ 0 $   &         &       \\
$ 1/2 $ & $ 1/2 $ &       \\ \hline
        & $ 0 $   & $ 1 $ 
\end{tabular} 
\end{center}
Since it is an order 2 method, for the computation of 
$ \operatorname{dexp}^{-1} $ one can use
$ \operatorname{dexp}^{-1}_\Omega (A) \approx A $.

\subsection{Numerical methods for Lie systems} \label{metodos}

So far, we have established in Procedure \ref{Metodo7Pasos} how to construct an analytical solution of a Lie system on a manifold $N$ via a Lie group action on  $N$, which is obtained by means of the integration of the VG Lie algebra of the Lie system. On the other hand, in Section \ref{LieGroupsMethods} we have reviewed some methods in the literature providing a numerical approximation of the solution of \eqref{eq.grupo} remaining in the Lie group $G$ (which accounts for their most remarkable geometrical property).

Now, let us explain how we combine these two elements to construct our new numerical methods, so we retrieve the solution of \eqref{Lie-Scheffer} on $N$. 
Let $ \varphi $  be the Lie group action  \eqref{def_accion} and consider the solution of the system \eqref{eq.grupo} such that $ Y(0) = I$. This solution permits us to retrieve the solution on $N$  of \eqref{Lie-Scheffer} for small values of $ t $, i.e., when a solution $ Y(t) $ of \eqref{eq.grupo} stays close to the neutral element and hence  the Lie group action $\varphi$ is properly defined. Numerically, we have shown that the solutions of \eqref{eq.grupo} can be provided through the approximations of \eqref{eq.algebraS}, say $\{\Omega_k\}_{k=0,\ldots,\mathcal{N}}$, and \eqref{GeneralLieGroupMethod}, as long as we stay close enough to the origin. As particular examples, we have picked the Magnus and RKMK methods in order to get $\{\Omega_k\}_{k=0,\ldots,\mathcal{N}}$ and, furthermore, the sequence $\{Y_k\}_{k=0,\ldots,\mathcal{N}}$. Next, we establish the scheme providing the numerical solution to Lie systems.

\begin{definition}\label{MetodosLS}
Let us consider a Lie system evolving on a manifold $N$ of the form 
\[
\frac{dx}{dt}=\sum_{\alpha=1}^r b_{\alpha}(t)X_{\alpha}(x),\quad x(a)=x_0,
\]
and let 
$$\frac{dY}{dt}=A(t)Y, \qquad A(t)=\sum_{\alpha=1}^rb_{\alpha}(t)M_{\alpha},
$$
be its associated automorphic Lie system. We define the numerical solution to the Lie system, i.e., $\{x_k\}_{k=0,\ldots,\mathcal{N}}$, via the algorithm given next.
\newpage
\begin{algorithm} {\rm Lie systems method}
\label{ForcedAlgorithm}
\begin{algorithmic}[1]
\State {\bf Initial data}: $\mathcal{N},\, h,\, x_0,\, A(t),\, Y_0=I ,\, \Omega_0={\bf 0} .$
 \State {\bf Numerically solve } 
 $
 \frac{d\Omega}{dt}= \mbox{dexp}^{-1}_{\Omega}A(t)
 $
 \State{\bf Output} $\{\Omega_k\}_{k=1,\ldots,\mathcal{N}}$
    \For {$k= 1,\ldots, \mathcal{N}-1$} 
    
\[
\begin{split}
Y_{k+1}&=e^{\Omega_k}Y_k,\\
x_{k+1}&=\varphi(Y_{k+1},x_k),
\end{split}
\]
    \EndFor
    \State  {\bf Output:} $(x_1,x_2,...,x_\mathcal{N}).$
\end{algorithmic}
  \end{algorithm}

\end{definition}

At this point, we would like to highlight an interesting geometric feature of this method. On the one hand, the discretization is based on the numerical solution of the automorphic Lie system underlying the Lie system, which, itself, is founded upon the geometric structure of the latter. This numerical solution remains on $G$, i.e., $Y_k\in G$ for all $k$, due to the particular design of the Lie group methods (as long as $h$ is small). Given this, our construction respects as well the geometrical structure of the Lie system, since, in principle, it evolves on a manifold $N$. We observe that the iteration 
\[
x_{k+1}=\varphi(Y_{k+1},x_k)
\]
leads to this preservation, since $x_{k+1}\in N$ as long as $Y_{k+1}\in G$ and $x_k\in N$ (we recall that $\varphi:G\times N\rightarrow N$). Note as well that the direct application of a one-step method \eqref{OneStepMethod} on a general Lie system \eqref{Lie-Scheffer} would destroy this structure. 

For future reference, in regards of the Lie group methods \eqref{GeneralLieGroupMethod}, we shall refer to \eqref{met.mag2} as Magnus 2, to \eqref{met.mag4} as Magnus 4 and to \eqref{RKMK4} as, simply, RKMK (we recall that this methods is order 4 convergent).

\section{Application to \texorpdfstring{SL$ (n, \mathbb{R}) $}{SL(n, R)}} \label{capSLN}
\subsection{\texorpdfstring{SL($ 2, \mathbb{R} $)}{SL(2, R)} and the  Riccati equation} \label{Ric1}
Let us recall the first-order Riccati equation over the real line $ \mathbb{R} $. One can check a comprehensive description of all the physical applications of this equation in \cite{araujo20}. The Riccati equation reads
\begin{equation}\label{Eq:Ricc}
\frac{dx}{dt} = b_0(t) + b_1(t) x + b_2(t) x^2 \, ,
\end{equation}
where $ b_0(t), b_1(t), b_2(t) $ are arbitrary $t$-dependent functions. The associated $t$-dependent vector field is
$X = b_0(t) X_0 + b_1(t) X_1 + b_2(t) X_2$, 
where
\begin{equation*}
X_0 = \partial_x , \qquad X_1 = x \partial_x , \qquad X_2 = x^2 \partial_x 
\end{equation*} 
and whose commutators are
\begin{equation} \label{comsl2}
[X_0, X_1] = X_0 , \qquad [X_0, X_2] = 2 X_1 , \qquad [X_1, X_2] = X_2 .
\end{equation} 
This proves that the Riccati equation is a Lie system related to a VG Lie algebra isomorphic to $\mathfrak{sl}(2,\mathbb{R})$. Thus, we employ the 7-step method \ref{Metodo7Pasos} to study its solutions. We choose the   basis $ \mathcal{V} = \{M_0, M_1, M_2 \} $ of $\mathfrak{sl}(2,\mathbb{R})$ to integrate the VG Lie algebra to a Lie group action of SL$(2,\mathbb{R})$ on $\mathbb{R}$. In more detail, 
\begin{equation*}
M_0 = \begin{pmatrix} 0 & 1 \\ 0 & 0 \end{pmatrix} , \qquad 
M_1 = \frac{1}{2} \begin{pmatrix} 1 & 0 \\ 0 & -1 \end{pmatrix} , \qquad 
M_2 = \begin{pmatrix} 0 & 0 \\ -1 & 0 \end{pmatrix} .
\end{equation*} 
Note that 
$$
[M_0,M_1]=-M_1,\qquad [M_0,M_2]=-2M_1,\qquad [M_1,M_2]=-M_2.
$$
We obtain the flows for the vector fields $ X_0, X_1 $ and $ X_2 $ by integrating them in terms of the real parameters  $ \lambda_0, \lambda_1, \lambda_2 $, respectively. Indeed, the flows of the vector fields $X_0,X_1,X_2$ read
\begin{equation*}
\Phi_0(\lambda_0, x_0) = \lambda_0 + x_0 \, , \qquad \Phi_1(\lambda_1, x_0)  = x_0 e^{\lambda_1}, \qquad \Phi_2(\lambda_2, x_0) = \frac{x_0}{1 - \lambda_2 x_0},
\end{equation*}
correspondingly. Using canonical coordinates of the second-kind, we can write  $ Y \in \rm{SL}(2, \mathbb{R}) $ near the neutral element as
\begin{equation}\label{Eq:Dec}
 Y = \exp(\lambda_0 M_0) \exp(\lambda_1 M_1) \exp(\lambda_2 M_2)\, . 
\end{equation}
We define the Lie group action $ \varphi : \rm{SL}(2, \mathbb{R}) \times \mathbb{R} \to \mathbb{R} $ through the equations
\[ \varphi(\exp(\lambda_i M_i), x) = \Phi_i(\lambda_i, x) \qquad i = 0, 1, 2. \]
Calculating the three exponential expressions in (\ref{Eq:Dec}) and comparing the expression with an arbitrary element $ Y \in \mathrm{SL}(2, \mathbb{R}) $ with parameters $\alpha\delta-\beta\gamma=1$, we have
\begin{equation*} 
Y = \begin{pmatrix}
\alpha & \beta \\
\gamma & \delta
\end{pmatrix} =
\begin{pmatrix}
e^{\lambda_1 / 2} - \lambda_0 \lambda_2 e^{-\lambda_1 / 2} & \lambda_0 e^{-\lambda_1 / 2} \\
-\lambda_2 e^{-\lambda_1 / 2}                              & e^{-\lambda_1 / 2}
\end{pmatrix} ,
\end{equation*}
from where the parameters $(\lambda_0,\lambda_1,\lambda_2)$ read
\begin{equation} \label{paramRic2} 
\lambda_0 = \frac{\beta}{\delta} , \qquad \lambda_1 = -2 \log \delta \, , \qquad \lambda_2 = -\frac{\gamma}{\delta} . 
\end{equation}
The action is obtained as 
\begin{equation*}
\varphi(Y, x_0) = \varphi(\exp(\lambda_0 M_0)\cdot \exp(\lambda_1 M_1) \cdot \exp(\lambda_2 M_2), x_0) = \Phi_0(\lambda_0, \Phi_1(\lambda_1, \Phi_2(\lambda_2, x_0))), 
\end{equation*}
and substituting the flows, 
\[ \varphi(Y, x_0) = \lambda_0 + \frac{x_0}{1 - \lambda_2 x_0} e^{\lambda_1} .\]
Now, substituting the parameters \eqref{paramRic2} and bearing in mind that for any $ Y \in \rm{SL}(2, \mathbb{R}) $ it is fulfilled that $ \alpha \delta - \gamma \beta = 1 $, we can reach the expression of the action that results in a homography \cite{hartshorne67}
\begin{equation}\label{HomoSL2}
    \varphi(Y, x_0) = \frac{\alpha x_0 + \beta}{\gamma x_0 + \delta} . 
\end{equation}

\subsubsection*{Exact solution}

It is interesting to note that if the $t$-dependent coefficients of the Lie system are constants, the matrix $Y$ associated with the linear system on the Lie group is $t$-independent and the solution of the automorphic Lie system can be easily retrieved.

For example, consider the Riccati equation with constant coefficients
\[ \frac{dx}{dt} = 1 + 2 x + x^2, \]
obtained by assuming $ b_0(t)= 1 $, $ b_1(t) = 2 $ and $ b_2(t) = 1 $ in (\ref{Eq:Ricc}). The system on the group \eqref{sistemaGM} associated with this Riccati equation reads 
\begin{equation} \label{coefCons}
\frac{dY}{dt} = A  Y, \qquad Y(t) = \begin{pmatrix} y_{11}(t) & y_{12}(t) \\ y_{21}(t) & y_{22}(t) \end{pmatrix} \in \mathrm{SL}(2,\mathbb{R}), \qquad Y(0) = I_2 , 
\end{equation}
where $ I_2 $ is the identity $2\times 2$ matrix  and  $ A(t) $ is
\[ A = \sum_{i = 0}^2 b_i M_i = \begin{pmatrix} 1 & 1 \\ -1 & -1 \end{pmatrix} . \]

If we write \eqref{coefCons} in the canonical form
\begin{equation*}
\frac{dy_{11}}{dt} = y_{11} + y_{21} , \quad 
\frac{dy_{12}}{dt} = y_{12} + y_{22} , \quad
\frac{dy_{21}}{dt} = -y_{11} - y_{21} , \quad
\frac{dy_{22}}{dt} = -y_{12} - y_{22} ,
\end{equation*}
or equivalently, $ d\boldsymbol{y} / dt = \Sigma \boldsymbol{y} $, where
\begin{equation*}
\boldsymbol{y} = \begin{pmatrix}
y_{11} \\
y_{12} \\
y_{21} \\
y_{22}
\end{pmatrix}, \qquad \Sigma = \begin{pmatrix}
1  & 0  &  1 & 0 \\
0  & 1  &  0 & 1 \\
-1 & 0  & -1 & 0 \\
0  & -1 & 0  & -1 
\end{pmatrix} , \qquad
\boldsymbol{y}(0) = \begin{pmatrix}
1 \\ 0 \\ 0 \\ 1
\end{pmatrix} .
\end{equation*}

The solution of the system reads
\[ \boldsymbol{y}(t) = \exp\left( \int_{0}^{t} \Sigma(\tau) d \tau \right) \boldsymbol{y}(0) = \exp\left( t \Sigma \right) \boldsymbol{y}(0) . \]
Observe that the matrix $ \Sigma $ is constant, so the integration is trivial. Also, since $\Sigma$ is nilpotent, the exponential is simply truncated at order 2. In this way, we obtain the solution:
\begin{equation*}
\boldsymbol{y}(t) = (t + 1, t, -t, 1 - t)^\mathsf{T} \qquad \Rightarrow \qquad 
Y(t) = \begin{pmatrix}
t + 1 & t \\
-t    & 1 - t
\end{pmatrix} .
\end{equation*}
Applying the Lie group action, we retrieve the solution of the original system:
\begin{equation*}
x(t) = \varphi\left( Y(t), x_0 \right) = \frac{t x_0 + x_0 + t}{1 - t - t x_0} , 
\end{equation*}
where $ x_0 $ is the initial condition.

\subsubsection*{Numerical example}

Let us now put into practice the numerical methods proposed in Definition \ref{MetodosLS}. For this matter, we consider  
\begin{equation} \label{ejemRiccati}
\frac{dx}{dt} = 2 t - \frac{x}{t} + \frac{x^2}{t^3} \, , \qquad t \geq 1.
\end{equation}
This is another  Riccati equation with $t$-dependent coefficients $b_0(t)=2t$, $b_1(t)=-1/t$ and $b_2(t)=1/t^3.$
Its solution is
\begin{equation}\label{ExSL2}
x(t) = \frac{2 t^3 - 2 t^2}{2 t - 1}  
\end{equation}
for the initial condition $ x(1) = 0 $.

\begin{figure}[h]
\begin{center}
\begin{tikzpicture}
\begin{groupplot}[group style = {group size = 2 by 1, horizontal sep = 40 pt},  width = 0.48\textwidth, height = 7 cm]
\nextgroupplot[ylabel = {$ x(t) $}, legend pos = north west, xlabel = $ t $, xmin = 1, xmax = 10]
	\addplot[domain = 1:10, line width = 1 pt]{(2 * x^3 - 2 * x^2)/(2 * x - 1)};            \addlegendentry{exact};
	\addplot[mark = square] table[x = t, y = mag2, col sep = semicolon]{datos/riccati.csv}; \addlegendentry{Magnus2};
	\addplot[mark = o     ] table[x = t, y = mag4, col sep = semicolon]{datos/riccati.csv}; \addlegendentry{Magnus4};
	\addplot[mark = otimes] table[x = t, y = RKMK, col sep = semicolon]{datos/riccati.csv}; \addlegendentry{RKMK};
	\node[] at (axis cs: 3.2, 80) {$ [h = 3] $};
\nextgroupplot[legend pos = north west, xlabel = $ t $, xmin = 1, xmax = 10]
	\addplot[domain = 1:10, line width = 1 pt]{(2 * x^3 - 2 * x^2)/(2 * x - 1)};              \addlegendentry{exact};
	\addplot[mark = square] table[x = t2, y = 2mag2, col sep = semicolon]{datos/riccati.csv}; \addlegendentry{Magnus2};
	\node[] at (axis cs: 3.2, 70) {$ [h = 0.5] $};
\end{groupplot}
\end{tikzpicture}
\caption{Exact vs. numerical solutions of \eqref{ejemRiccati} with $ x(1) = 0 $. In the left plot we observe the natural better approximation of higher-order  methods for huge time steps ($h=3$). In the right plot, we observe a closer approximation to the exact dynamics when $h$ decreases ($h=0.5$). }
\label{ricc1}
\end{center}
\end{figure}
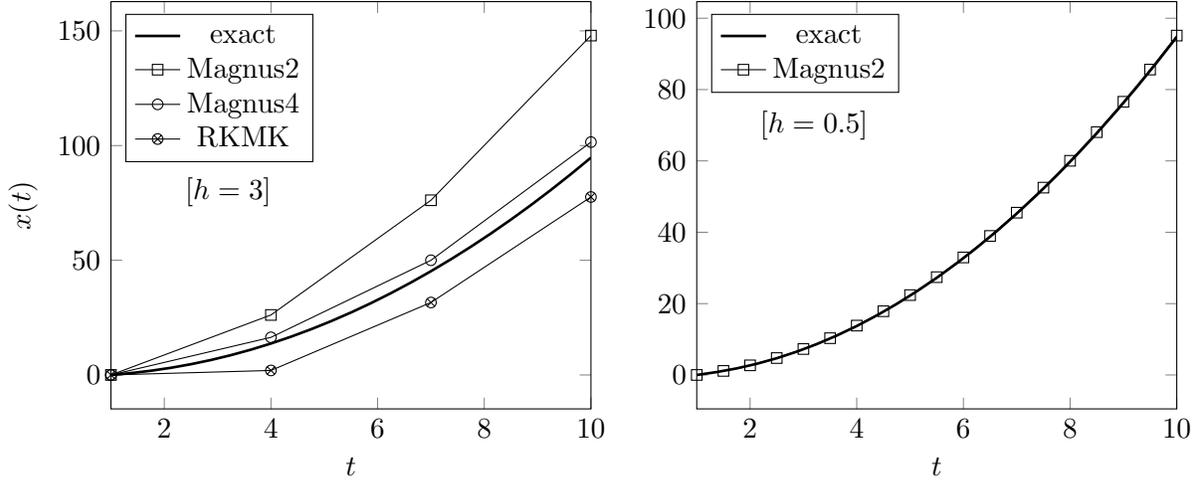

In Figure \ref{ricc1} we show how the described numerical methods approximate the exact solution \eqref{ExSL2} in the interval $ [1, 10] $ taking different time steps and employing Magnus 2, Magnus 4 and RKMK as underlying methods in the Lie group.

\begin{figure}[h]
\begin{center}
\begin{tikzpicture}
\begin{groupplot}[group style = {group size = 2 by 1, horizontal sep = 40 pt},  width = 0.48\textwidth]
\nextgroupplot[ylabel = {$ E_N $}, xmode = log, ymode = log, legend pos = north west]
	\addplot[mark = o     ] table[x = h2, y = mag2, col sep = semicolon]{datos/convRiccati.csv}; \addlegendentry{Magnus2};
	\addplot[mark = square] table[x = h2, y = Heun, col sep = semicolon]{datos/convRiccati.csv}; \addlegendentry{Heun};   
\nextgroupplot[xmode = log, ymode = log, legend pos = north west]
	\addplot[mark = square] table[x = h4, y = RKMK, col sep = semicolon]{datos/convRiccati.csv}; \addlegendentry{RKMK};
	\addplot[mark = o     ] table[x = h4, y = mag4, col sep = semicolon]{datos/convRiccati.csv}; \addlegendentry{Magnus4};
	\addplot[mark = otimes] table[x = h4, y = RK4,  col sep = semicolon]{datos/convRiccati.csv}; \addlegendentry{RK4};
\end{groupplot}
\end{tikzpicture}
\caption{Convergence for the Riccati equation \eqref{ejemRiccati}.}
\label{convRic}
\end{center}
\end{figure}
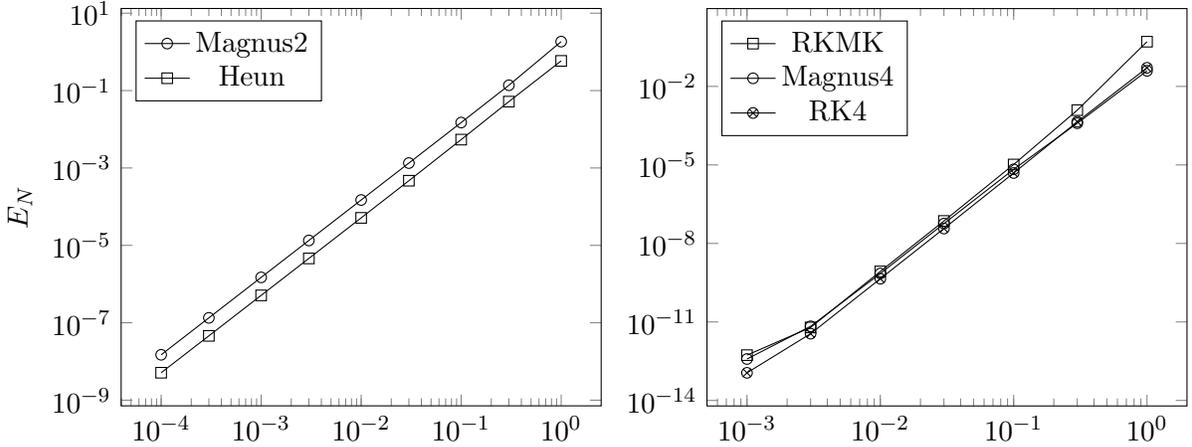

In Figure  \ref{convRic} we show convergence plots. To make a proper comparison we include two classical numerical schemes, Heun (order 2) and RK4, respectively, for the corresponding orders, applied directly to \eqref{ejemRiccati}.  As it is apparent, the slope of the convergence lines are two and four, and this manifests that the order of convergence of the numerical methods on the underlying Lie group  is transmitted to the manifold in this particular example. This transmission can be easily understood in terms of the local truncation error of the underlying Lie group method and the particular form of the analytical solution we obtain, i.e., \eqref{HomoSL2}. Namely, if we are applying an order $p$ Lie group method in this particular example, that means $\alpha_{k+1}=\alpha(t_{k+1})+\mathcal{O}(h^{p_1+1})$, $\beta_{k+1}=\beta(t_{k+1})+\mathcal{O}(h^{p_2+1})$, $\gamma_{k+1}=\gamma(t_{k+1})+\mathcal{O}(h^{p_3+1})$, $\delta_{k+1}=\delta(t_{k+1})+\mathcal{O}(h^{p_4+1})$, where $p=$min$\{ p_1,\,p_2,\,p_3,\,p_4\}$. Naturally, $\alpha,\,\beta,\,\gamma,\,\delta$ are the components of the SL$(2,\mathbb{R})$ matrix we are dealing with. Taking this into account, the analytical expression \eqref{HomoSL2} and the definition of the local truncation error we have introduced in \eqref{LTE}, it is straightforward to see that $E_h=\mathcal{O}(h^{p+1})$, and, consequently, it is to expect that the convergence order of the Lie group method is transmitted to the manifold.

\subsection{\texorpdfstring{SL($ 3, \mathbb{R} $)}{SL(3, R)} and matrix Riccati equations} \label{riccatiM}

A general matrix Riccati equation \cite{HWA83} has the following form
\begin{equation}\label{matrixric}
 \frac{d \Gamma}{dt} = G_1(t) + G_2(t) \Gamma + \Gamma G_3(t) + \Gamma G_4(t) \Gamma,
 \end{equation}
where 
$\Gamma, G_1(t) \in \mathcal{M}_{n \times m}(\mathbb{R}),\,G_2(t) \in \mathcal{M}_{n \times n}(\mathbb{R}),\, G_3(t) \in \mathcal{M}_{m \times m}(\mathbb{R}),\, G_4(t) \in \mathcal{M}_{m \times n}(\mathbb{R}).$ The case that matters to us is
 $ n = 2, m = 1 $, for which the matrix Riccati equation has a VG Lie algebra isomorphic $\mathfrak{sl}(3,\mathbb{R})$. Then, equation \eqref{matrixric} takes the form
\begin{equation}\label{matrixRicc}
\begin{pmatrix}
\dot{x} \\ \dot{y}
\end{pmatrix} = 
\begin{pmatrix}
g_1(t) \\ g_2(t)
\end{pmatrix} + 
\begin{pmatrix}
g_3(t) & g_4(t) \\
g_5(t) & g_6(t)
\end{pmatrix}
\begin{pmatrix}
x \\ y
\end{pmatrix} + 
\begin{pmatrix}
x \\ y
\end{pmatrix} g_7(t) + 
\begin{pmatrix}
x \\ y
\end{pmatrix}
\begin{pmatrix}
g_8(t) & g_9(t)
\end{pmatrix}
\begin{pmatrix}
x \\ y
\end{pmatrix},
\end{equation}
where $ g_1(t),\ldots,g_9(t)$ are arbitrary functions of time. Equivalently, we can write the previous matrix equation as
\begin{equation}\label{Eq:OurMRE}
\left\{
\begin{aligned}
\frac{dx}{dt} &= g_1(t) + \left( g_3(t) + g_7(t) \right) x + g_4(t) y + g_8(t) x^2 + g_9(t) xy, \\
\frac{dy}{dt} &= g_2(t) + g_5(t) x + \left( g_6(t) + g_7(t) \right) y + g_8(t) xy +  g_9(t) y^2.
\end{aligned}
\right.
\end{equation}
The $t$-dependent vector field associated with this system can be written as
\begin{equation*}
X = g_1(t)X_1 + g_2(t) X_2 + \left( g_3(t) + g_7(t) \right) X_3 + \left( g_6(t) + g_7(t) \right) X_4 + g_4(t) X_5 + g_5(t) X_6 + g_8(t) X_7 + g_9(t) X_8 ,
\end{equation*}
where
\begin{align*}
X_1 = \partial_x, \quad \quad X_2 = \partial_y, \quad & \quad X_3 = x \partial_x, \quad  \quad X_4 = y \partial_y, \\
X_5 = y \partial_x, \quad  \quad X_6 = x \partial_y, \quad  \quad X_7 &= x^2 \partial_x + xy \partial_y, \quad  \quad X_8 = xy \partial_x + y^2 \partial_y.
\end{align*}
Note that   $ X$ only really depends on eight $t$-dependent functions, since $ g_3(t), g_6(t) $ and $ g_7(t) $ appear as linear combinations $ g_3 (t)+ g_7(t) $  and $ g_6(t) + g_7(t) $. Let us list only the non-vanishing commutators for these vector fields:
\begin{equation}\label{commsl3}
\begin{gathered}[]
[X_1, X_3] = X_1, \quad [X_1, X_6] = X_2, \quad [X_1, X_7] = 2 X_3 + X_4, \quad [X_1, X_8] = X_5, \\
[X_2, X_4] = X_2, \quad [X_2, X_5] = X_1, \quad [X_2, X_7] = X_6, \quad [X_2, X_8] = X_3 + 2 X_4, \\
[X_3, X_5] = -X_5, \quad [X_3, X_6] = X_6, \quad [X_3, X_7] = X_7, \quad [X_4, X_5] = X_5, \\
[X_4, X_6] = -X_6, \quad [X_4, X_8] = X_8, \quad [X_5, X_6] = X_4 - X_3, \\
[X_5, X_7] = X_8, \quad [X_6, X_8] = X_7. 
\end{gathered}
\end{equation}

From this we conclude that \eqref{matrixRicc} is a Lie system. Now, we choose a matrix basis for $ \mathfrak{sl}(3, \mathbb{R}) $:
\begin{equation*}
\begin{gathered}
M_1 =             \begin{pmatrix}  0 & 0 & 0 \\ 0 &  0 & 0 \\ -1 &  0 &  0 \end{pmatrix} , \quad
M_2 =             \begin{pmatrix}  0 & 0 & 0 \\ 0 &  0 & 0 \\  0 & -1 &  0 \end{pmatrix} , \quad  
M_3 = \frac{1}{3} \begin{pmatrix}  2 & 0 & 0 \\ 0 & -1 & 0 \\  0 &  0 & -1 \end{pmatrix} , \quad
M_4 = \frac{1}{3} \begin{pmatrix} -1 & 0 & 0 \\ 0 &  2 & 0 \\  0 &  0 & -1 \end{pmatrix} , \\[3 pt]
M_5 =             \begin{pmatrix}  0 & 0 & 0 \\ 1 &  0 & 0 \\  0 &  0 &  0 \end{pmatrix} , \quad  
M_6 =             \begin{pmatrix}  0 & 1 & 0 \\ 0 &  0 & 0 \\  0 &  0 &  0 \end{pmatrix} , \quad
M_7 =             \begin{pmatrix}  0 & 0 & 1 \\ 0 &  0 & 0 \\  0 &  0 &  0 \end{pmatrix} , \quad
M_8 =             \begin{pmatrix}  0 & 0 & 0 \\ 0 &  0 & 1 \\  0 &  0 &  0 \end{pmatrix} .
\end{gathered}
\end{equation*}
To integrate the VG Lie algebra of (\ref{Eq:OurMRE}) to a Lie group action, we express the elements of the Lie group $ \mathrm{SL}(3, \mathbb{R}) $ in terms of canonical coordinates of the second-kind in the following way
\begin{equation}\label{Eq:SL3SCC}
Y = \begin{pmatrix} y_{11} & y_{12} & y_{13} \\ y_{21} & y_{22} & y_{23} \\ y_{31} & y_{32} & y_{33} \end{pmatrix} = \prod_{i = 1}^8 \exp(\lambda_i M_i)\in \mathrm{SL}(3, \mathbb{R}),
\end{equation} 
where $ \lambda_1, \dots, \lambda_8 \in \mathbb{R} $ are real parameters univocally determined for each $Y$ in an open neighborhood of the neutral element of ${\rm SL}(3,\mathbb{R})$. The exponentials in the above expression can be calculated very easily and by using their values in (\ref{Eq:SL3SCC}) it turns our that
\begin{equation}\label{Movidon1}
\begin{split}
y_{11} = k_1 , \quad y_{12}& = \lambda_6 k_1 , \quad y_{13} = (\lambda_6 \lambda_8 + \lambda_7) k_1 , \\ 
y_{21} = \lambda_5 k_2 , \quad y_{22} = (1 + \lambda_5 \lambda_6) k_2 , &\quad y_{23} = (\lambda_5 \lambda_7 + \lambda_5 \lambda_6 \lambda_8 + \lambda_8) k_2 , \\
y_{31} = -\lambda_1 k_1 - \lambda_2 \lambda_5 k_2 , &\quad y_{32} = -\lambda_2 (1 + \lambda_5 \lambda_6) k_2 - \lambda_1 \lambda_6 k_1 , \\
y_{33} = -\lambda_1 (\lambda_6 \lambda_8 + \lambda_7) k_1 - \lambda_2 (\lambda_5 \lambda_7& + \lambda_5 \lambda_6 \lambda_8 + \lambda_8) k_2 + e^{(-\lambda_3 -\lambda_4) / 3} ,
\end{split}
\end{equation}
where $ k_1 = e^{(2 \lambda_3 - \lambda_4) / 3} $ y $ k_2 = e^{(2 \lambda_4 - \lambda_3) / 3} $. Rewriting some equalities in terms of others, i.e., $
y_{31} = -\lambda_1 y_{11} - \lambda_2 y_{21}$ and $ y_{32} = -\lambda_1 y_{12} - \lambda_2 y_{22}$,
we obtain a linear system from where we get $ \lambda_1 $ and $ \lambda_2 $. Operating with the remaining ones, we calculate the rest of the parameters. 
\begin{equation}\label{Movidon2}
\begin{split}
\lambda_1     = \frac{y_{22} y_{31} - y_{21} y_{32}}{y_{12} y_{21} - y_{11} y_{22}} , &\qquad
\lambda_2     = \frac{y_{11} y_{32} - y_{12} y_{31}}{y_{12} y_{21} - y_{11} y_{22}} , \\[5 pt]
e^{\lambda_3} = -y_{11}(y_{12} y_{21} - y_{11} y_{22}) , &\qquad 
e^{\lambda_4} = \frac{(y_{12} y_{21} - y_{11} y_{22})^2}{y_{11}} , \\[5 pt]
\lambda_5 = -\frac{y_{11} y_{21}}{y_{12} y_{21} - y_{11} y_{22}} , &\qquad
\lambda_6 = \frac{y_{12}}{y_{11}} , \\[5 pt]
\lambda_7 = \frac{y_{12} y_{23} - y_{13} y_{22}}{y_{12} y_{21} - y_{11} y_{22}} , &\qquad
\lambda_8 = \frac{y_{13} y_{21} - y_{11} y_{23}}{y_{12} y_{21} - y_{11} y_{22}} .
\end{split}
\end{equation} 
Integrating the vector fields $X_1,\ldots,X_8$, we obtain their flows, $\Phi_1,\ldots,\Phi_8$, which in turn give us the action 
\begin{gather*}
\begin{aligned}
\Phi_1(\lambda_1, (x_0, y_0)) = (\lambda_1 + x_0, y_0) , &\qquad
\Phi_2(\lambda_2, (x_0, y_0)) = (x_0, \lambda_2 + y_0) , \\[3 pt]
\Phi_3(\lambda_3, (x_0, y_0)) = (x_0 e^{\lambda_3}, y_0) , &\qquad
\Phi_4(\lambda_4, (x_0, y_0)) = (x_0, y_0 e^{\lambda_4}) , \\[3 pt]
\Phi_5(\lambda_5, (x_0, y_0)) = (x_0 + y_0 \lambda_5, y_0) , &\qquad
\Phi_6(\lambda_6, (x_0, y_0)) = (x_0, y_0 + x_0 \lambda_6) , \\[3 pt]
\end{aligned} \\
\Phi_7( \lambda_7 , (x_0, y_0)) = \left( \frac{x_0}{1 - x_0 \lambda_7}, \frac{y_0}{1 - x_0 \lambda_7} \right) , \quad
\Phi_8( \lambda_8 , (x_0, y_0)) = \left( \frac{x_0}{1 - y_0 \lambda_8}, \frac{y_0}{1 - y_0 \lambda_8} \right) .
\end{gather*}
In view of (\ref{Eq:SL3SCC}),  the composition of the flows $ \Phi_1 \circ \Phi_2 \circ \cdots \circ \Phi_8 $ allows us to obtain the complete action $(x,y)= \varphi(Y, (x_0, y_0))$, with
\begin{equation*}
x = \frac{x_0 (1 + \lambda_5 \lambda_6) + y_0 \lambda_5}{1 - x_0 \lambda_7 - y_0 \lambda_8} e^{\lambda_3} + \lambda_1 , \qquad 
y = \frac{x_0 \lambda_6 + y_0}{1 - x_0 \lambda_7 - y_0 \lambda_8} e^{\lambda_4} + \lambda_2 .
\end{equation*}

Operating with these expressions, we can rewrite the action through homographies as follows
\begin{equation*}
(x_0, y_0) \to \left( \frac{a_{21} + a_{22} x_0 + a_{23} y_0}{a_{11} + a_{12} x_0 + a_{13} y_0}, \frac{a_{31} + a_{32} x_0 + a_{33} y_0}{a_{11} + a_{12} x_0 + a_{13} y_0} \right) , 
\end{equation*}
with coefficients
\begin{equation}\label{Movidon3}
\begin{gathered}
a_{11} = 1, \qquad  \qquad a_{12} = -\lambda_7, \qquad  \qquad a_{13} = -\lambda_8 ,\qquad  \qquad a_{21} = \lambda_1, \\
a_{22} = (1 + \lambda_5 \lambda_6) e^{\lambda_3} - \lambda_1 \lambda_7, \qquad  \qquad a_{23} = \lambda_5 e^{\lambda_3} - \lambda_1 \lambda_8, \qquad \qquad a_{31} = \lambda_2, \\
a_{32} = \lambda_6 e^{\lambda_4} - \lambda_2 \lambda_7, \qquad \qquad a_{33} = e^{\lambda_4} - \lambda_2 \lambda_8 .
\end{gathered}
\end{equation}

\subsubsection*{Numerical example}
To illustrate  again our numerical methods we will take the following equation as an example:
\begin{equation} \label{ejemRiccM}
\left\lbrace \begin{aligned}
\frac{dx}{dt} &= 5 \sin 10 t - x + y , \\
\frac{dy}{dt} &= 5 \cos 10 t + x + y ,
\end{aligned} \right. 
\end{equation}
which is a matrix Riccati equation \eqref{matrixRicc} with $t$-dependent functions
\begin{gather*}
g_1(t) = 5 \sin 10 t, \qquad g_2(t) = 5 \cos 10 t, \qquad g_3(t)+ g_7(t) = -1 , \\
g_4(t) = 1, \qquad g_5 = 1, \qquad g_6(t) + g_7(t) = 1, \qquad g_8(t) = g_9(t) = 0 .
\end{gather*}
More exactly, it is an affine system of first-order differential equations. For the initial condition $ (x(0), y(0)) = (1, 0) $, the solution of \eqref{ejemRiccM} is
\begin{equation*} 
\left\lbrace \begin{aligned}
x(t) &= \frac{157}{102} \cosh \sqrt{2} t - \frac{\sqrt{2}^3 19}{51} \sinh \sqrt{2} t + \frac{5}{102} \sin 10 t - \frac{55}{102} \cos 10 t , \\
y(t) &= \frac{27 \sqrt{2}}{34} \sinh \sqrt{2} t + \frac{5}{102} \cosh \sqrt{2} t + \frac{15}{34} \sin 10 t - \frac{5}{102} \cos 10 t .
\end{aligned} \right. 
\end{equation*}

Figure 3 shows convergence plots.

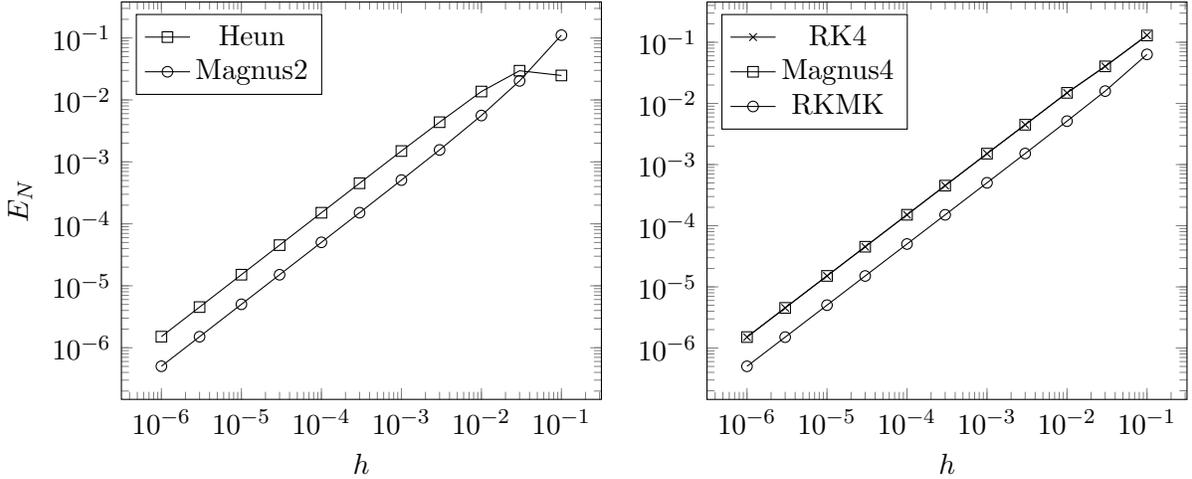
\begin{figure}[ht]
\begin{center}
\begin{tikzpicture}
\begin{groupplot}[group style = {group size = 2 by 1, horizontal sep = 40 pt},  width = 0.48\textwidth]
\nextgroupplot[xlabel = $ h $, ylabel = $ E_N $, xmode = log, ymode = log, legend pos = north west]
	\addplot[mark = square] table[x = h2, y = Heun, col sep = semicolon]{datos/convRiccatiM.csv}; \addlegendentry{Heun};
	\addplot[mark = o     ] table[x = h2, y = mag2, col sep = semicolon]{datos/convRiccatiM.csv}; \addlegendentry{Magnus2};
\nextgroupplot[xlabel = $ h $, xmode = log, ymode = log, legend pos = north west]
	\addplot[mark = x] table[x = h4, y = RK4,  col sep = semicolon]{datos/convRiccatiM.csv};      \addlegendentry{RK4};
	\addplot[mark = square] table[x = h4, y = mag4, col sep = semicolon]{datos/convRiccatiM.csv}; \addlegendentry{Magnus4};
	\addplot[mark = o     ] table[x = h4, y = RKMK, col sep = semicolon]{datos/convRiccatiM.csv}; \addlegendentry{RKMK};
\end{groupplot}
\end{tikzpicture}
\caption{Convergence for the affine system of first-order differential equations \eqref{ejemRiccM}.}
\label{convRiccM}
\end{center}
\end{figure}

In this case one can depict that, although our method is still convergent, the order of the Lie group is not transmitted to the manifold $N$ (in both cases the slope of the convergence lines is about  1). In this case, our method is not compared to Heun and RK4 applied directly to \eqref{ejemRiccM} but to an alternate scheme given by
\begin{equation}\label{Alternativa}
x_{k+1}=\varphi(\tilde Y_{k+1},x_k),
\end{equation}
where $\{ \tilde Y_k\}_{k=0,\ldots,\mathcal{N}}$ is the numerical solution of \eqref{eq.grupo} when Heun and RK4 are applied to them (in Figure 3 they are referred as Heun and RK4). Naturally, this implies that $\tilde Y_k\notin G$.  

Our conjecture is that, in this case, the construction of the action changes the convergence of the method, which can be sustained in the high nonlinearity obtained when defining the parameters \eqref{Movidon1}, \eqref{Movidon2}, \eqref{Movidon3}. An interesting open question is whether there is a way to modify the methods according to the Lie group action so that the convergence is transmitted correctly. Another clue pointing in that direction is that, as it can be easily seen in the plots, although the velocity of convergence is about the same for our method and \eqref{Alternativa}, quantitatively the error of the former is lower. We consider this as another (positive) geometrical symptom, since, apparently, the error  worsens when the underlying Lie group structure is not preserved.



\subsection{Generalization to \texorpdfstring{SL($ n, \mathbb{R} $)}{SL(n, R)}} \label{genera}

The special linear Lie group plays an essential role in mechanical systems and integrable systems (see \cite{hussin90,LG17,Re72} and references therein). This is why we briefly detail a possible generalization of our proposed methods to SL($ n, \mathbb{R} $).

Recall that the Lie algebra $ \mathfrak{sl}(n, \mathbb{R}) $ associated with the Lie group $ \mathrm{SL}(n, \mathbb{R}) $ has dimension $ n^2 - 1 $. In fact,
a matrix representation of $\mathfrak{sl}(n,\mathbb{R})$ is given by the matrix Lie algebra given by $ n \times n $ traceless matrices. For simplicity, we can choose a basis of $\mathfrak{sl}(n,\mathbb{R})$ given by 
 $ n^2 - n $ matrices with one nontrivial off-diagonal entry equal to one, together with  $ n - 1 $ diagonal traceless matrices of the form
\begin{equation*} 
\begin{pmatrix}
1      &  0     & 0      & \cdots & 0      \\
0      & -1     & 0      & \cdots & 0      \\
0      &  0     & 0      & \cdots & 0      \\
\vdots & \vdots & \vdots & \ddots & \vdots \\
0      & 0      & 0      & \cdots & 0 
\end{pmatrix} , \quad
\begin{pmatrix}
1      &  0     & 0      & \cdots & 0      \\
0      &  0     & 0      & \cdots & 0      \\
0      &  0     & -1     & \cdots & 0      \\
\vdots & \vdots & \vdots & \ddots & \vdots \\
0      & 0      & 0      & \cdots & 0 
\end{pmatrix} ,\quad \cdots \quad
\begin{pmatrix}
1      &  0     & 0      & \cdots & 0      \\
0      &  0     & 0      & \cdots & 0      \\
0      &  0     & 0      & \cdots & 0      \\
\vdots & \vdots & \vdots & \ddots & \vdots \\
0      & 0      & 0      & \cdots & -1
\end{pmatrix} .
\end{equation*}

The total $ n^2 - 1 $ matrices are traceless and linearly independent. A Lie group action  $\varphi:\mathrm{SL}(n, \mathbb{R})\times \mathbb{R}^{n-1}\rightarrow \mathbb{R}^{n-1}$ can be then constructed via homographies as follows (cf. \cite{HWA83}):

\begin{equation} \label{homografia}
\begin{aligned}
x_1       &= \frac{a_{10} + a_{11} x_1^0 + a_{12} x_2^0 + \cdots a_{1, n - 1} x_{n - 1}^0}{a_{00} + a_{01} x_1^0 + a_{02} x_2^0 + \cdots a_{0, n - 1} x_{n - 1}^0} , \\[7 pt]
x_2       &= \frac{a_{20} + a_{21} x_1^0 + a_{22} x_2^0 + \cdots a_{2, n - 1} x_{n - 1}^0}{a_{00} + a_{01} x_1^0 + a_{02} x_2^0 + \cdots a_{0, n - 1} x_{n - 1}^0} , \\
          & \hspace{6 pt} \vdots \\
x_{n - 1} &= \frac{a_{n - 1, 0} + a_{n - 1, 1} x_1^0 + a_{n - 1, 2} x_2^0 + \cdots a_{n - 1, n - 1} x_{n - 1}^0}{a_{00} + a_{01} x_1^0 + a_{02} x_2^0 + \cdots a_{0, n - 1} x_{n - 1}^0} ,
\end{aligned}
\end{equation}
where $ (x_1, \dots, x_{n - 1})=\varphi(Y,(x_1^0, \dots, x_{n - 1}^0) )$, where  $ Y \in \mathrm{SL}(n, \mathbb{R}) $ is
\begin{equation*} 
Y = \begin{pmatrix}
a_{00}       & a_{01}       & a_{02}       & \cdots & a_{0, n - 1}    \\
a_{10}       & a_{11}       & a_{12}       & \cdots & a_{1, n - 1}    \\
a_{20}       & a_{21}       & a_{22}       & \cdots & a_{2, n - 1}    \\
\vdots       & \vdots       & \vdots       & \ddots & \vdots          \\
a_{n - 1, 0} & a_{n - 1, 1} & a_{n - 1, 2} & \cdots & a_{n - 1, n - 1} 
\end{pmatrix} .
\end{equation*}
Note that if $ \langle \, \cdot \, , \cdot \, \rangle $ is the standard scalar product in $\mathbb{R}^n$ and we call $ a_i $, with $i=0,\ldots,n$, the rows of $ Y $ and  $ \overline{x}^0 $ stands for  the point $ (1, x_1^0, \dots, x_{n - 1}^0) $ in $ \mathbb{R}^n $, then \eqref{homografia} can be rewritten as $ x_i = \langle a_i, \overline{x}^0 \rangle / \langle a_0, \overline{x}^0 \rangle $ for $ i = 1, \dots, n - 1 $. 

It is worth noting that if two VG Lie algebras $V_1,V_2$ on two manifolds $N_1,N_2$ are diffeomorphic, i.e. there exists a diffeomorphism $\phi:N_1\rightarrow N_2$ such that $\phi_*V_1=V_2$, then $V_1,V_2$ can be integrated to two $\phi$-equivariant Lie group actions $\varphi_1:G\times N_1\rightarrow N_1$ and $\varphi_2:G\times N_2\rightarrow N_2$, i.e., $\phi(\varphi_1(g,x))=\varphi_2(g,\phi(x))$ for every $x\in N_1$ and $g\in G$. In particular if $V_1$ is the VG Lie algebra of matrix Riccati equations studied in this section and $V_2$ is another VG Lie algebra on $N_2=N_1$ diffeomorphic to $V_1$, then the Lie group action $\varphi_2$ is $\phi$-equivariant to $\varphi_1$. Since every diffeomorphism in $N_1$ can be understood as a change of variables, the $\phi$-equivariance of $\varphi_1$ and $\varphi_2$ entails that a change of variables in $N_2$ allows us to write the action of every $g\in {\rm SL}(n,\mathbb{R})$ via $\varphi_2$ as an homography. Note that it is simple to prove that \eqref{homografia} gives rise to a Lie group action of SL$(n,\mathbb{R})$ and its fundamental vector fields are those related to matrix Riccati equations. 

\subsubsection{Increase of numerical cost as $n$ increases}

We can indirectly measure the numerical cost of our schemes according to the time they need to compute the solution. Let us consider the following equation
\begin{equation*} 
	\frac{dx}{dt} = 2 t - \frac{x}{t} + \frac{x^2}{t^3} \, , \qquad t \geq 1 , \qquad x(1) = 0 ,
\end{equation*}
whose analytical solution is
\[ x(t) = \frac{2 t^3 - 2 t^2}{2 t - 1} \, . \]

 Now, we apply our five numerical schemes to the equation above and plot the step size (which is strictly related with $ n $) versus the time consumed for the resolution of the equation. 
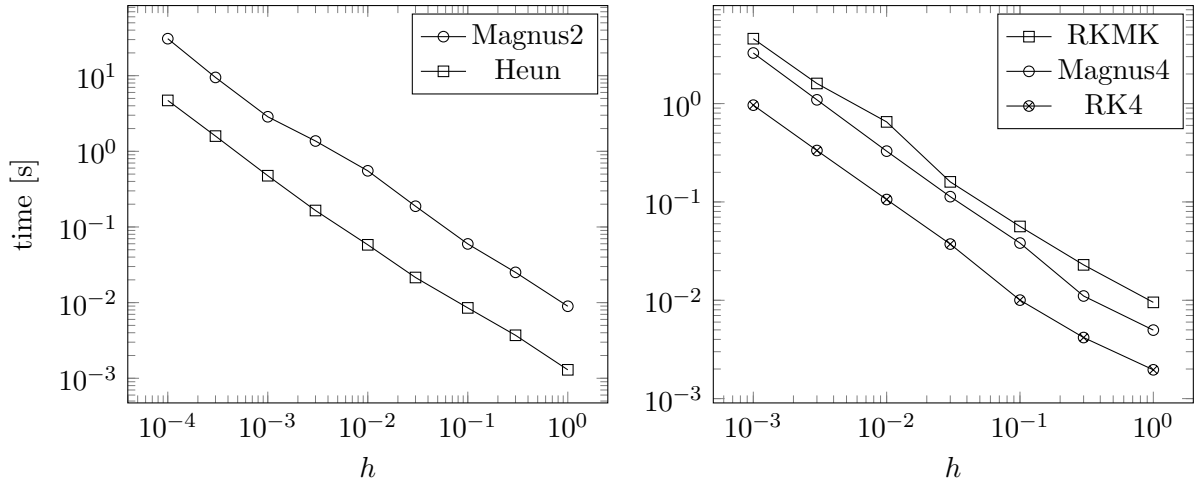
\begin{center}
	\begin{tikzpicture}
		\begin{groupplot}[group style = {group size = 2 by 1, horizontal sep = 40 pt},  width = 0.48\textwidth]
			\nextgroupplot[xlabel = $ h $, ylabel = {time [s]}, xmode = log, ymode = log]
			\addplot[mark = o     ] table[x = h2, y = mag2, col sep = semicolon]{tiempoRiccati.csv}; \addlegendentry{Magnus2};
			\addplot[mark = square] table[x = h2, y = Heun, col sep = semicolon]{tiempoRiccati.csv}; \addlegendentry{Heun};
			\nextgroupplot[xlabel = $ h $, xmode = log, ymode = log]
			\addplot[mark = square] table[x = h4, y = RKMK, col sep = semicolon]{tiempoRiccati.csv}; \addlegendentry{RKMK};
			\addplot[mark = o     ] table[x = h4, y = mag4, col sep = semicolon]{tiempoRiccati.csv}; \addlegendentry{Magnus4};
			\addplot[mark = otimes] table[x = h4, y = RK4, col sep = semicolon]{tiempoRiccati.csv}; \addlegendentry{RK4};
		\end{groupplot}
	\end{tikzpicture}
\end{center}
We can observe that in the logarithmic axis the relation between the variables is close to being linear. As expected, the 4th-order schemes (RKMK, Magnus 4 and RK4) show a bigger increase in numerical cost as $ n $ increases. 

Now, we renact the same process to the following differential system
\begin{equation*} 
	\left\lbrace \begin{aligned}
		\frac{dx}{dt} &= 5 \sin 10 t - x + y , \\
		\frac{dy}{dt} &= 5 \cos 10 t + x + y ,
	\end{aligned} \right. \qquad 
	\left\lbrace \begin{aligned}
	x(0) &= 1, \\
	y(0) &= 0,
\end{aligned} \right.
\end{equation*}
whose solution can be written as
\begin{equation*} 
	\left\lbrace \begin{aligned}
		x(t) &= \frac{157}{102} \cosh \sqrt{2} t - \frac{\sqrt{2}^3 19}{51} \sinh \sqrt{2} t + \frac{5}{102} \sin 10 t - \frac{55}{102} \cos 10 t , \\
		y(t) &= \frac{27 \sqrt{2}}{34} \sinh \sqrt{2} t + \frac{5}{102} \cosh \sqrt{2} t + \frac{15}{34} \sin 10 t - \frac{5}{102} \cos 10 t .
	\end{aligned} \right. 
\end{equation*}
and we obtain the following graphics. 
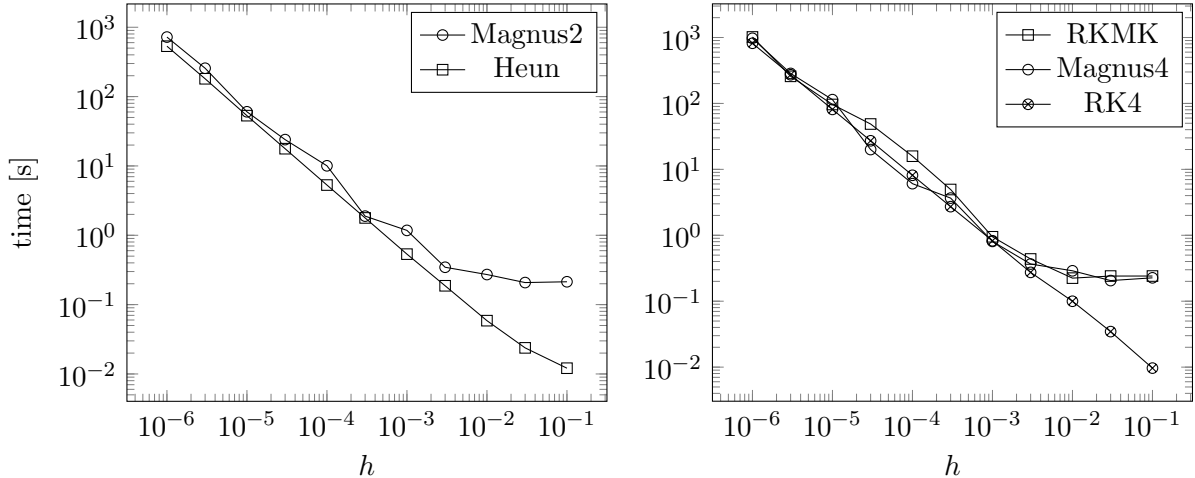
\begin{center}
	\begin{tikzpicture}
		\begin{groupplot}[group style = {group size = 2 by 1, horizontal sep = 40 pt},  width = 0.48\textwidth]
			\nextgroupplot[xlabel = $ h $, ylabel = {time [s]}, xmode = log, ymode = log]
			\addplot[mark = o     ] table[x = h2, y = mag2, col sep = semicolon]{tiempoRiccatiM.csv}; \addlegendentry{Magnus2};
			\addplot[mark = square] table[x = h2, y = Heun, col sep = semicolon]{tiempoRiccatiM.csv}; \addlegendentry{Heun};
			\nextgroupplot[xlabel = $ h $, xmode = log, ymode = log]
			\addplot[mark = square] table[x = h4, y = RKMK, col sep = semicolon]{tiempoRiccatiM.csv}; \addlegendentry{RKMK};
			\addplot[mark = o     ] table[x = h4, y = mag4, col sep = semicolon]{tiempoRiccatiM.csv}; \addlegendentry{Magnus4};
			\addplot[mark = otimes] table[x = h4, y = RK4,  col sep = semicolon]{tiempoRiccatiM.csv}; \addlegendentry{RK4};
		\end{groupplot}
	\end{tikzpicture}
	\label{tiempoRiccM}
\end{center}
When $ h $ is small (and, therefore $ n $ is big) we observe again a linear relation between the numerical cost and the index $ n $.  

\section{Applications in Linear Quadratic Control} \label{control}

Now, we provide an interesting application to optimal control of the method to obtain the solution of Lie systems given in Procedure \ref{Metodo7Pasos}. A very useful model to carry out the control of dynamical systems is the representation in the space of states. The most general representation on such space is
\begin{equation*}
\left\lbrace \begin{aligned}
\dot{\boldsymbol{x}}(t) &= \boldsymbol{f}(\boldsymbol{x}(t), \boldsymbol{u}(t), t) , \\
\boldsymbol{y}(t)       &= \boldsymbol{h}(\boldsymbol{x}(t), \boldsymbol{u}(t), t) ,
\end{aligned} \right. 
\end{equation*} 
where $ \boldsymbol{x} : [t_0, t_f] \to \mathbb{R}^n $ is a vector containing the state variables of the system, $ \dot{\boldsymbol{x}} : [t_0, t_f] \to \mathbb{R}^n $ is its time derivative, $ \boldsymbol{u} : [t_0, t_f] \to \mathbb{R}^m $ is the vector containing the input variables, $ \boldsymbol{y} : [t_0, t_f] \to \mathbb{R}^p $ is the vector with the output variables and $ \boldsymbol{f} : \mathbb{R}^n \times \mathbb{R}^m \times \mathbb{R} \to \mathbb{R}^n $ and $ \boldsymbol{h} : \mathbb{R}^n \times \mathbb{R}^m \times \mathbb{R} \to \mathbb{R}^p $ are two $t$-dependent arbitrary vector fields. We can manipulate the inputs to modify the state of the system. 

A very important and common model is that of linear systems, given their simplicity \cite{dominguez06}. Indeed, it is pretty usual to search for a linearization of nonlinear problems. The most general representation of a linear system is
\begin{equation} \label{sistLineal}
\left\lbrace  \begin{aligned}
\dot{\boldsymbol{x}}(t) &= A(t) \boldsymbol{x}(t) + B(t) \boldsymbol{u}(t) , \\
\boldsymbol{y}(t)       &= C(t) \boldsymbol{x}(t) + D(t) \boldsymbol{u}(t) ,
\end{aligned} \right. 
\end{equation} 
where the $t$-dependent matrices $ A(t), B(t), C(t) $ and $ D(t) $ are the state (or system) matrix, the input matrix, the output matrix and the feedthrough (or feedforward) matrix, respectively. In order for the system to be defined the dimensions of the matrices must be $ A(t) \in \mathcal{M}_n $, $ B(t) \in \mathcal{M}_{n \times m} $, $ C(t) \in \mathcal{M}_{p \times n} $ and $ D(t) \in \mathcal{M}_{p \times m} $ for every $t\in \mathbb{R}$.

In particular, we are interested in the problem of optimal control with a quadratic cost function, which, as we are going to show, can be transformed into a matrix Ricatti equation. This is, given a linear system
\eqref{sistLineal}, the state $ \boldsymbol{x}_0 $ and the time interval $ [t_0, t_f] $, we need to find an input $ \boldsymbol{u}(t) $ starting with condition $ \boldsymbol{x}(t_0) = \boldsymbol{x}_0 $ that minimizes the quadratic cost function, i.e.,
\begin{equation*} 
J(\boldsymbol{x}, \boldsymbol{u}) \overset{\text{def.}}{=} \boldsymbol{x}(t_f)^\textsf{T} S  \boldsymbol{x}(t_f) + \int_{t_0}^{t_f} \boldsymbol{x}(t)^\textsf{T} Q(t) \boldsymbol{x}(t) dt + \int_{t_0}^{t_f} \boldsymbol{u}(t)^\textsf{T} R(t) \boldsymbol{u}(t) dt ,
\end{equation*} 
where $ S  $ is a positive semi-definite matrix and for all $ t \in [t_0, t_f] $ the matrices  $ Q(t) $ and $ R(t) $ are, respectively, positive semi-definite and positive definite. Obviously, $ S , Q(t) \in \mathcal{M}_{n\times n} $ and $ R(t) \in \mathcal{M}_{m\times m} $ for every $t\in \mathbb{R}$.

Since the matrices involved are positive (semi-) definite, the terms appearing in them are a measure of the size of the  vectors $ \boldsymbol{x} $ and $ \boldsymbol{u} $. Each of them ``penalizes" a different aspect of the control.
The first one measures how far the system is from the null state $ \boldsymbol{x} = 0 $ at the end of the time interval.
Analogously, the second term measures the distance between the state and the null state along time. In this way, the fastest the system approaches the null state, the smallest the cost function is and the closest it is to the null state at the end of the time interval.
On the other hand, the third term measures the size of the input along time in such a way that the smallest it is (with respect to the measure defined by the matrix $R(t)$), the smallest the value of the function $J$ will be.

Adjusting the matrices $ S , Q(t), R(t) $ we choose what aspects are more important. If we choose the matrix $ S $ in such a way that staying far from the null state at the end of the interval is very penalized, the optimal control will conduct the system towards this state at the end of the time interval, at the cost that the input will be bigger. If $ Q(t) $ takes over the other two matrices, the control will lead the system to the null state as fast as possible. On the contrary, if the dominant matrix is $ R(t) $, the input $ \boldsymbol{u} $ will be small, but probably the other two aspects will be adversely affected. This is interesting when the size of the input is related to any other variable that we would like to minimize. 

In this formulation the cost function leads the system towards the null state. Nonetheless, it is easy to modify the problem so the system drifts towards a different state.
If we aim at establishing the system in a certain state $ \boldsymbol{x}_c $, if we are capable of finding an input $ \boldsymbol{u}_c $ such that
\[ 0 = A \boldsymbol{x}_c + B \boldsymbol{u}_c , \]
then, performing the change of variables
\begin{align*} \left\lbrace \begin{aligned}
\boldsymbol{x}(t) - \boldsymbol{x}_c &\to \hat{\boldsymbol{x}}(t), \\
\boldsymbol{u}(t) - \boldsymbol{u}_c &\to \hat{\boldsymbol{u}}(t) ,
\end{aligned} \right. 
\end{align*} 
we obtain a new system
\begin{equation*}
\frac{d \hat{\boldsymbol{x}}}{dt} = \frac{d \boldsymbol{x}}{dt} = A(\hat{\boldsymbol{x}} + \boldsymbol{x}_c) + B(\hat{\boldsymbol{u}} + \boldsymbol{u}_c) = 
A \hat{\boldsymbol{x}} + B \hat{\boldsymbol{u}} + \underbrace{A \boldsymbol{x}_c + B \boldsymbol{u}_c}_{= 0} = A \hat{\boldsymbol{x}} + B \hat{\boldsymbol{u}}
\end{equation*}
in which we can apply the quadratic cost function to obtain an optimal control problem that conducts the system towards
$ \hat{\boldsymbol{x}} = \boldsymbol{x} - \boldsymbol{x}_c = 0 $. In this way, the original system will tend to $ \boldsymbol{x}_c $.

The solution of the linear quadratic control problem is given as a state-feedback controller, i.e., the optimal input $ \boldsymbol{u}_o(t) $ that minimizes $ J(\boldsymbol{x}, \boldsymbol{u}) $ is a function of the state of the system.  In particular, we can write $ \boldsymbol{u}_o(t) = K(t) \boldsymbol{x}(t) $, where $ K(t) $ is the feedback matrix and it is calculated as 
\[ K(t) = -R(t)^{-1} B(t)^\textsf{T} P(t) , \]
where $ P(t) $ is the solution of the following matrix differential Riccati equation
\begin{equation} \label{controlOpt}
\frac{d P}{dt} = P B R(t)^{-1} B^\textsf{T} P - P A - A^\textsf{T} P - Q(t), \qquad P(t_f) = S .
\end{equation}
The initial condition is given at the end of the time interval because one needs to integrate the equation in reverse \cite[\S 8.2]{sontag98}. Equation \eqref{controlOpt} is the matrix Riccati equation introduced in Section \ref{capSLN}. 

Now, we are going to solve an example involving linear quadratic control by the application of our analytical resolution of Lie systems.

\subsection{Example: velocity of a vehicle}

We propose a  model of a control for the velocity of a vehicle. We will have a single input variable, which will correspond with the strength of the engine to accelerate the vehicle. Let us assume that the only force that could decelerate the vehicle is the friction with air and that it is proportional to the square of the velocity \cite{PGG17}. For simplicity, our model reduces to describing motions with positive velocity. Under these hypotheses, it is enough to take the velocity of the vehicle as the variable of state to completely characterize the system. Applying the second law of Newton, we obtain the equation describing the system
\[ F - k v^2 = m \frac{dv}{dt} , \]
where $ F $ is the engine force, $ v $ is the velocity, $ k $ is a constant of proportionality and $ m $ is the mass of the vehicle. For simplicity, we will take   $ m = k = 1 $. We change the notation to use $ u $ instead of $ F $, being this one the input of the system. So, the system reads now
\[ \frac{dv}{dt} = -v^2 + u . \]
This system is nonlinear, but when we are designing a control that keeps the velocity constant around a certain value, we can linearize the system in a neighborhood of such value to compute the optimal control with quadratic cost function that keeps the vehicle at cruising speed. Again, to simplify the computations we take $ v_c = 1 $.
Under these circumstances, $ dv/dt = 0 $ so we obtain $ u_c = 1 $. The linearized system around the point $ (v_c, u_c) $ results in
\[ \frac{d \Delta v}{dt} = -2 \Delta v + \Delta u , \]
where $ \Delta v = v - 1 $ and $ \Delta u = u - 1 $ are the incremental variables around $ (v_c, u_c) $.

 To further simplify, we will take all the matrices constant in the quadratic cost function, and equal to one in the time interval $ [0, 1] $. Then, the cost function is
\begin{equation}\label{costfunctionexample}
 J(v, u) = \Delta v(1)^2 + \int_0^1 \Delta v(s)^2 ds + \int_0^1 \Delta u(s)^2 ds . 
\end{equation}
The function $ \Delta u_o(t) $ that minimizes \eqref{costfunctionexample} is  $ \Delta u_o = K(t) \Delta v(t) $, where $ K = -R^{-1} B^\textsf{T} P $. In our case $ K(t) = -P(t) $, being $ P(t) $ the solution of the Riccati equation
\begin{equation} \label{ejemLQ}
\frac{d P}{dt} = P B R^{-1} B^\textsf{T} P - P A - A^\textsf{T} P - Q = P^2 + 4 P - 1 , 
\end{equation}
with (final) condition $ P(1) = S = 1 $.

Now, we resolve \eqref{ejemLQ} analytically, applying our procedure exposed in \ref{Metodo7Pasos}. Since it is a Riccati equation with constant coefficients, given its simplicity, we can compute its analytical solution by resolving its associated the linear system  on the group $ \mathrm{SL}(2, \mathbb{R}) $. In this case, we have to solve $ dY(t)/dt = A Y(t) $, with $ Y(1) = I $, where the matrix $A$ is (according to the notation in Section \ref{Ric1})
\[ A = -M_0 + 4 M_1 + M_2 = \begin{pmatrix} 2 & -1 \\ -1 & -2 \end{pmatrix} . \] 
The exact solution of this system will be expressed in its canonical form $ d\boldsymbol{y}/dt = \Sigma \boldsymbol{y} $, where
\[ \Sigma = \begin{pmatrix} 2 & 0 & -1 & 0 \\ 0 & 2 & 0 & -1 \\ -1 & 0 & -2 & 0 \\ 0 & -1 & 0 & -2 \end{pmatrix}, \qquad \boldsymbol{y}(1) = \begin{pmatrix} 1 \\ 0 \\ 0 \\ 1 \end{pmatrix} . \]
Its solution is
\[ \boldsymbol{y}(t) = \exp \left( \int_{1}^{t} \Sigma d \tau \right) \boldsymbol{y}(1) = \exp \left( (t - 1) \Sigma \right) \boldsymbol{y}(1), \]
this is,
\begin{equation*}
\boldsymbol{y}(t) = \frac{e^{-(\sqrt{5} t + \sqrt{5})}}{10}
\begin{pmatrix} 
    (5 + 2 \sqrt{5}) e^{2 \sqrt{5} t} + (5 - 2 \sqrt{5}) e^{2 \sqrt{5}} \\[5 pt]
    \sqrt{5} e^{2 \sqrt{5}} - \sqrt{5} e^{2 \sqrt{5} t} \\[5 pt]
    \sqrt{5} e^{2 \sqrt{5}} - \sqrt{5} e^{2 \sqrt{5} t} \\[5 pt]
    (5 - 2 \sqrt{5}) e^{2 \sqrt{5} t} + (5 + 2 \sqrt{5}) e^{2 \sqrt{5}} 
\end{pmatrix} .
\end{equation*}

Finally, we can retrieve the solution to \eqref{ejemLQ} by means of the Lie group action of $ \mathrm{SL}(2, \mathbb{R}) $ on $ \mathbb{R} $ as
\[ P(t) = \varphi\left(Y(t), P(1) \right) = \frac{(5 + \sqrt{5}) e^{2 \sqrt{5} t} + (5 - \sqrt{5}) e^{2 \sqrt{5}}}{(5 - 3 \sqrt{5}) e^{2 \sqrt{5} t} + (5 + 3 \sqrt{5}) e^{2 \sqrt{5}}} . \]

The optimal control problem is $ \Delta u_o = -P(t) \Delta v $. We introduce a constant $ \Delta u_c $ that carries the system from an initial perturbation to the functioning point $ v = 1 $.
If we start from a point $ v(0) = \overline{v} $, to determine the constant value of $ u $ that takes the system back to the cruising speed $ \Delta u_c $ the equation
\[ \frac{d \Delta v}{dt} = -2 \Delta v + \Delta u_c,  \]
with initial conditions $ \Delta v(0) = \overline{v} - 1 $ and $ \Delta v(1) = 0 $. The solution can be computed trivially
\[ \Delta u_c = \frac{2 \overline{v} - 2}{1 - e^2} . \]

\begin{figure}
	\begin{center}
		\begin{tikzpicture}
		\begin{axis}[xmin = 0, xmax = 1, width = 0.95\textwidth, height = 5.6 cm, legend entries={optimal, constant}, xlabel = $ t $, ylabel = $ v(t) $]
		\addlegendimage{};
		\addlegendimage{gray, dashed};
		\addplot[domain = 0:1]{1};            
		\addplot[] table[x = tiempo, y = v1.2, col sep = semicolon]{datos/control.csv}; 
		\addplot[] table[x = tiempo, y = v0.9, col sep = semicolon]{datos/control.csv}; 
		\addplot[gray, dashed] table[x = tiempo, y = c1.2, col sep = semicolon]{datos/control.csv}; 
		\addplot[gray, dashed] table[x = tiempo, y = c0.9, col sep = semicolon]{datos/control.csv};
		\end{axis}
		\end{tikzpicture}
		\caption{Evolution of the system with optimal control and constant control}
		\label{figCont}
	\end{center}
\end{figure}
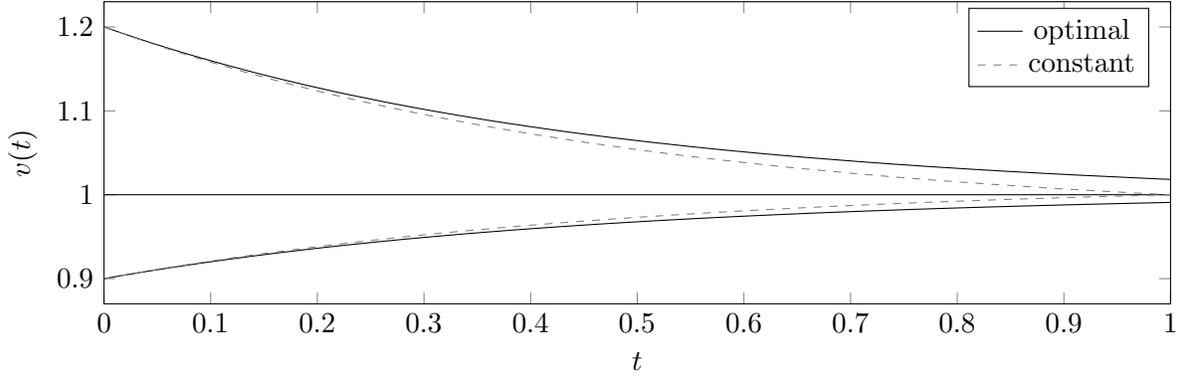

In Figure \ref{figCont} we have depicted the evolution of the system with different initial conditions around $ v = 1 $.
The continuous line represents the evolution of the system when we use optimal control and the discontinuous line corresponds with a constant $ u $.  

The chosen values of $ S, Q, R $ for the optimal control do not take the vehicle at cruising speed in the time interval considered. This makes sense if we think of the quadratic cost function as a compromise to reduce the size of the input, so the system reaches the functioning point fast and efficiently. If we want to ensure that the vehicle reaches the cruising speed, we need to reflect it in the cost function by giving more weight to  $ S $ and $ Q $.  

If we now calculate the cost function for different initial conditions, we see that the constant control makes the system reach the cruising speed quicker and with less error than the optimal control, and the cost is smaller. We list some values on the following table.

\begin{center}
\begin{tabular}{l|ccccc}
	    $ J(v, u) \times 1000 $ & $ \overline{v} = 1.2 $ & $ \overline{v} = 1.15 $ & $ \overline{v} = 1.1 $ & $ \overline{v} = 1.05 $ & $ \overline{v} = 1 $ \\ \hline
	    Optimal Control          & 10,340                 & 5,816                   & 2,585                   & 0,646                  & 0                    \\
	    Constant Control       & 11,771                 & 6,621                   & 2,943                   & 0,736                  & 0 
\end{tabular} 
\end{center}


The input represents the engine force  accelerating the vehicle. It is also reasonable that the fuel consumption will be proportional to the strength of the force.
In this way, we can derive the optimal control that keeps the vehicle at constant cruising speed and minimizing the amount of fuel. 

\section{Conclusions}

This paper is concerned with the integration of Lie systems, both from the analytical and numerical perspectives,  using particular techniques adapted to their geometric features.
This work is rooted in the field of numerical and discrete methods specifically adapted for Lie systems which is still a very unexplored brach of research \cite{araujo20,PW04,piet12,RW84}. 

One major result in this paper is that we are able to solve Lie systems on Lie groups. This permits us to solve all Lie systems related to the same automorphic Lie system at the same time (equivalently, all Lie systems that have isomorphic VG Lie algebras) \cite{carinena00,lucas11}. Automorphic Lie systems present a simple superposition rule that only depends on a single particular solution. This is an advantage in comparison with superposition rules for general Lie systems, which use to depend on a larger number of particular solutions. The second most important advantage is that, since  Lie groups admit a local matrix representation, automorphic Lie systems can be written as first-order systems of linear homogeneous ODEs in normal form.

{Employing the geometric structure of Lie systems, we propose a particular geometric integrator for Lie systems that exploits the properties of such structure}. Particularly, we employ the Lie group action obtained by integrating the Vessiot-Guldberg Lie algebra of a Lie system to get the analytical solution of the Lie system. We use the automorphic Lie system related to a Lie system, along with geometric schemes, say Lie group integrators, to preserve the group structure. Specifically, we use two families of numerical schemes: the first one based on the Magnus expansion, whereas the second is based on RKMK methods. We have compared both methods in different situations. We can generally say that the fourth-order RKMK is slightly more precise than the Magnus expansion of the same order. Regarding the transmission of convergence order from the Lie group method to the Lie system method, our conjecture, rooted in the results obtained for different Lie groups, is that how the Lie group action is constructed has a central role. Whilst the numerical methods work very satisfactorily on the Lie group level, when we translate the properties into the manifold we see that the convergence and precision of the numerical method can be modified ({as in the SL$(3,\mathbb{R})$ case}). Nonetheless, since our methods are based on geometric integrators, they inherit all the geometric properties we wish to preserve and the solutions always belong in the manifold, where the Lie system is defined (something that is not preserved if one uses classical numerical schemes).

From the results obtained for $ \mathrm{SL}(2, \mathbb{R}) $ and $ \mathrm{SL}(3, \mathbb{R}) $ we have been able to provide a generalization to $ \mathrm{SL}(n, \mathbb{R}) $, and we have discussed the form of the Lie group action. As it has been evidenced, $ \mathrm{SL}(n, \mathbb{R}) $ is a relevant Lie group, appearing recurrently in nonlinear oscillators of Winternitz-Smorodinsky, Milney-Pinney, Ermakov sytems, as well as higher-order Riccati equations.

The last important result is that solving higher-order Riccati equations has allowed us to resolve important examples appearing in engineering problems. We have particularly proposed a problem in optimal control in which matrix Riccati equations appear naturally from quadratic cost functions.

{In the future, we will analyse the convergence transmission from automorphic Lie systems to related Lie systems. In addition, since the exponential is a local diffeomorphism, the topological study of matrix Lie groups would allow us to establish the optimal time-step for Lie group methods, which is a long-standing problem that would also help optimize Lie system methods. Another endevour is to study Lie systems on more general manifolds that are not necessarily isomorphic to $\mathbb{R}^n$ and depict how some geometric and topological invariants are preserved \cite{HLT17,LL19}.}
Right now, we are working on examples on Anti-de-Sitter spaces so we can depict how the curvature is preserved under the numerical method. We could easily generalize this to all kinds of systems in all types of curved spaces.
This will in fact prove the interest of our 7-step method, since one could argue that the nongeometric approximation methods seem fairly better than our proposal. Nonetheless, in our forthcoming publications we will show that when there are invariants in the game, the 7-step method is the best choice to preserve certain geometric and topological invariants.
\section*{Acknowledgements}
 J. de Lucas   acknowledges partial financial support from  MINIATURA-5 Nr 2021/05/X/ST1/01797, funded by the National Science Centre (Poland). C. Sardón and F. Jiménez acknowledge project “Teoría de aproximación constructiva y aplicaciones" (TACA-ETSII), UPM, Madrid.
 
 \section*{Data availability}
 The datasets generated during and/or analysed during the current study are available from the corresponding author on reasonable request.


\begin{thebibliography}{60}

\bibitem{ado47} I.D. Ado, \textit{The representation of Lie algebras by matrices}. Uspekhi Matematicheskikh Nauk, \textbf{2}:159-173, 1947.

\bibitem{angelo05} R.M. Angelo and W.F. Wreszi\'nski, \textit{Two-level quantum dynamics, integrability and unitary NOT gates}, Phys. Rev. A \textbf{72}:034105, 2005.

\bibitem{blanes00} S. Blanes, F. Casas and J. Ros, \textit{Improved high order integrators based on the Magnus expansion}, BIT Nume. Math. \textbf{40}:434-450, 2000. 

\bibitem{blasco15} A. Blasco, F.J. Herranz, J. de Lucas and C. Sardon, \textit{Lie-Hamilton systems on the plane: Applications and superposition rules}, J. Phys. A  \textbf{48}:345202, 2015.

\bibitem{carinena00} J.F. Cariñena, J. Grabowski and G. Marmo, {\sl Lie-Scheffers Systems: a Geometric Approach}, Napoli Series in Physics and Astrophysics, Bibliopolis, 2000.

\bibitem{carinena07} J.F. Cariñena, J. Grabowski and G. Marmo, \textit{Superposition rules, Lie theorem and partial differential equations}, Rep. Math. Phys. \textbf{60}:237-258, 2007.

\bibitem{carinena01} J.F. Cariñena, J. Grabowski and A. Ramos, \textit{Reduction of $t$-dependent systems admitting a superposition principle}, Acta Appl. Math. \textbf{66}:67-87, 2001.

\bibitem{carinena09} J.F. Cariñena and J. de Lucas, \textit{Applications of Lie systems in dissipative Milne-Pinney equations}, Int. J. Geom. Meth.  Modern Phys. \textbf{6}:683-699, 2009.

\bibitem{lucas11} J.F. Cariñena and J. de Lucas, {\sl Lie Systems: Theory, Generalisations, and Applications}, Dissertationes Mathematicae \textbf{479}, 2011.

\bibitem{carinena11} J.F. Cariñena, J. de Lucas and C. Sardón, \textit{A new Lie systems approach to second-order Riccati equations}, Int. J.  Geom. Meth.   Modern Phys. \textbf{9}:1260007, 2011.

\bibitem{carinena99} J.F. Cariñena and A. Ramos, \textit{Integrability of the Riccati equation from a group theoretical viewpoint}, Int. J. Modern Phys. A \textbf{14}:1935-1951, 1999.


\bibitem{CortesMartinez} J. Cort\'es and S. Mart\'inez, \textit{Non-holonomic integrators}, Nonlinearity \textbf{14}:1365-1392, 2001.

\bibitem{MLC}
 M.L. Curtis, {\sl Matrix groups} 2nd ed. New York: Springer, cop. 1984. Universitext. 0387960740
\bibitem{dominguez06} S. Domínguez, P. Campoy, J.M. Sebastián and A. Jiménez,  {\sl Control en el Espacio de Estado}, Pearson, Educación, 2006.




\bibitem{hairer06} E. Hairer, C. Lubich and G. Wanner, {\sl Geometric Numerical Integration}, Springer-Verlag, Berlin-Heidelberg, 2006.

\bibitem{hairer93} E. Hairer, S.P. N\o rsett and G. Wanner, {\sl Solving Ordinary Differential Equations I: Nonstiff Problems}, Springer-Verlag, Berlin Heidelberg, 1993.

\bibitem{Hall}
 B. Hall, Matrix Lie Groups in {\sl Lie Groups, Lie Algebras, and Representations: An Elementary Introduction},3--30, Springer International Publishing, 2015.




\bibitem{HWA83}
J. Harnad, P. Winternitz and R.L. Anderson, 
{\it Superposition principles for matrix Riccati equations}, 
J. Math. Phys. {\bf 24}:1062, 1983.
\bibitem{hartshorne67} R. Hartshorne, {\sl Foundations of Projective Geometry}, W.A. Benjamin, Inc., Nueva York, 1967.


\bibitem{HLT17}
 F.J. Herranz, J. de Lucas and M. Tobolski, {\it 
Lie-Hamilton systems on curved spaces: A geometrical approach},
J. Phys. A {\bf 50}:495201, 2017. 

\bibitem{hussin90} V. Hussin, J. Beckers, L. Gagnon and P. Winternitz, \textit{Superposition formulas for nonlinear superequations}, J. Math. Phys. \textbf{31}:2528-2534, 1990.

\bibitem{isaacson66} E. Isaacson and H.B. Keller, {\sl Analysis of Numerical Methods}, John Wiley \& Sons, New York-London-Sydney, 1966.

\bibitem{iserles05} A. Iserles, H. Munthe-Kaas, S. N\o rsett and A. Zanna, {\it Lie-group methods}, Acta Numerica 215-365, 2005.

\bibitem{iserles99} A. Iserles and S.P. N\o rsett, \textit{On the solution of linear differential equations in Lie groups}, Phil. Trans Royal Soc. A \textbf{357}:983-1020, 1999.

\bibitem{iserles98} A. Iserles, S.P. N\o rsett and A.F. Rasmussen, \textit{$t$-symmetry and high-order Magnus methods}, Technical Report 1998/NA06, DAMTP, University of Cambridge, 1998.

\bibitem{Ku73}
V. Ku\v{c}era, 
{\it A Review of the Matrix Riccati Equation},
Kybernetika {\bf 9}:42-61, 1973.
\bibitem{LL19}
 J. Lange and J. de Lucas, {\it Geometric models for Lie--Hamilton systems on $\mathbb{R}^2$},
 Mathematics {\bf 2019}:7, 1053.
 
\bibitem{lazaro09} J.A. Lázaro-Camí and J.P. Ortega, \textit{Superposition rules and stochastic Lie-Scheffers systems}, Ann. Inst. H. Poincaré Probab. Stat. \textbf{45}:910-931, 2009.

\bibitem{lee00} J.M. Lee, {\sl Introduction to Smooth Manifolds}, Graduate Texts in Mathematics 218, Springer-Verlag, New York, 2003. 

\bibitem{levi05} E.E. Levi, \textit{Sulla struttura dei gruppi finiti e continui}, Atti della Reale Accademia delle Scienze di Torino, 1905.

\bibitem{lie93} S. Lie and G. Scheffers, {\sl Vorlesungen über continuierliche Gruppen mit geometrischen und anderen Anwendungen}, Teubner, Leipzig, 1893.

\bibitem{LG17}
J. de Lucas and A.M. Grundland,
{\it A Lie systems approach to the Riccati hierarchy and partial differential equations}, J.  Differential Equations {\bf 263}:299-337 (2017).


 
\bibitem{araujo20} J. de Lucas and C. Sardón, {\sl A Guide to Lie Systems with Compatible Geometric Structures}, World Scientific, Singapore, 2020. 



\bibitem{Magnus} W. Magnus, \textit{On the exponential solution of differential equations for a linear operator}, Comm. Pure Appl. Math. \textbf{7}:649-673, 1954.

\bibitem{MMM} J.C. Marrero, D. Martín de Diego and E. Mart\'inez, \textit{Discrete Lagrangian and Hamiltonian mechanics on Lie groupoids}, Nonlinearity \textbf{19}:1313-1348, 2006.


\bibitem{MaWe} J.E. Marsden and M. West, \textit{Discrete mechanics and variational integrators}, Acta Numerica \textbf{10}:357-514, 2001.


\bibitem{McLac} R. McLachlan and G.R.W. Quispel, \textit{Splitting methods}, Acta Numerica \textbf{11}:341-434, 2002.

\bibitem{munthe98} H. Munthe-Kaas, \textit{Runge-Kutta methods on Lie groups}, BIT Numerical Mathematics \textbf{38}:92-111, 1998.

\bibitem{munthe99} H. Munthe-Kaas, \textit{High order Runge-Kutta methods on manifolds}, J. Appl. Num. Maths. \textbf{29}:115-127, 1999.

\bibitem{odzijewicz00} A. Odzijewicz and A.M. Grundland, \textit{The Superposition Principle for the Lie Type first-order PDEs}, Rep. Math. Phys. \textbf{45}:293-306, 2000.
\bibitem{PGG17}
A. Pandey, A. Ghose-Choudhury and P. Guha,
{\it 
Chiellini integrability and quadratically damped oscillators},
Int. J. Non-Linear Mechanics {\bf 92}:153-159, 2017.

\bibitem{PW04}
A.V. Penskoi and P. Winternitz, {\it Discrete matrix Riccati equations with super- position formulas}, J. Math. Anal. Appl. {\bf 294}:533–547, 2004.


\bibitem{piet12} G. Pietrzkowski, \textit{Explicit solutions of the $ \mathfrak{a}_1 $-type Lie-Scheffers system and a general Riccati equation},  J. Dyn. Control Systems \textbf{18}:551-571, 2012.

\bibitem{quarteroni07} A. Quarteroni, R. Sacco and F. Saleri, {\sl Numerical Mathematics}, Springer-Verlag, New York, 2007.

\bibitem{RW84}
D.W. Rand and P. Winternitz, {\it Nonlinear superposition principles: a new numerical method for solving matrix Riccati equations}, Comput. Phys. Comm. {\bf 33}:305-328, 1984.


\bibitem{Re72} W.T. Reid, {\sl Riccati Differential Equations}, Academic, New York, 1972.

\bibitem{SS} J.M. Sanz-Serna, \textit{Symplectic integrators for Hamiltonian problems: an overview}, Acta Numerica {243-286}, 1992.

\bibitem{sardon15} C. Sardón, {\sl Lie systems, Lie symmetries and reciprocal transformations}, PhD Thesis, Universidad de Salamanca, 2015.

\bibitem{Sattinger}

D.H. Sattinger, O.L. Weaver, {\sl Lie Groups and Algebras with Applications to Physics}, Geometry and Mechanics. Berlin, Heidelberg: Springer-Verlag, 1986. Applied Mathematical Sciences, 61. 9781441930774

\bibitem{sontag98} E.D. Sontag, {\sl Mathematical Control Theory: Deterministic Finite Dimensional Systems}, Springer-Verlag, New York, 1998.

\bibitem{win82} P. Winternitz, \textit{Nonlinear action of Lie groups and superposition rules for nonlinear differential equations}, Phys. A {\bf 114}:105-113, 1982.


\bibitem{Zanna} A. Zanna, \textit{Collocation and relaxed collocation for the Fer and Magnus expansions}, J. Numer. Anal. \textbf{36}:1145-1182, 1999.

\end{thebibliography}
\end{document}